\documentclass[a4paper]{amsart}
\usepackage{fullpage}
\usepackage{amsmath,amssymb,mathtools,booktabs,tikz}
\usetikzlibrary{arrows.meta,calc,decorations.markings,positioning}
\usepackage[T1]{fontenc}
\usepackage{lmodern}
\usepackage{microtype}
\usepackage{needspace}
\usepackage{tabularx,array}
\usepackage[hypertexnames=false]{hyperref}
\usepackage{doi}
\usepackage{float}
\hypersetup{
    colorlinks=true,
    linkcolor=blue,
    citecolor=blue,
    urlcolor=blue
}
\usepackage{graphicx}
\newtheorem{thm}{Theorem}[section]
\newtheorem{lem}[thm]{Lemma}

\newtheorem{cor}[thm]{Corollary}
\newtheorem{conj}[thm]{Conjecture}

\theoremstyle{definition}
\newtheorem{rem}[thm]{Remark}
\numberwithin{equation}{section}
\allowdisplaybreaks[2]
\newcommand{\D}{\overline{\mathbb D}}
\newcommand{\eps}{\varepsilon}
\newcommand{\Err}{\mathcal E}
\newcommand{\dd}{\,\mathrm d}
\newcommand{\ee}{\mathrm e}

\DeclareMathOperator{\Rea}{Re}

\title[A proof of Schmeisser's conjecture]{A proof of Schmeisser's conjecture}
\author[T.~Zhang]{Teng Zhang}
\address{School of Mathematics and Statistics, Xi'an Jiaotong University, Xi'an 710049, P. R. China}
\email{teng.zhang@stu.xjtu.edu.cn}
\date{November 13, 2026}

\begin{document}

\begin{abstract}
In 1977, Schmeisser proposed a strengthening of Sendov's
conjecture: if all zeros of a polynomial $p$ of degree
$n\ge2$ lie in the closed unit disk, then for every point
$\xi$ in the convex hull of the zeros of $p$, there exists
a critical point $\zeta$ of $p$ such that $
    |\xi-\zeta|\le 1.$
We prove this conjecture in every degree.
The proof combines a covering of the convex hull by
diameter disks, inequalities for reflected polynomials,
and quantitative estimates for centred products,
supplemented by rigorous computer-assisted verification
in a finite range of degrees.  The Lean 4 formalization of the main results are also provided.
\end{abstract}

\maketitle

\section{Introduction}

In 1958, Sendov conjectured that if all zeros of a polynomial
$p$ of degree $n\ge2$ lie in the closed unit disk, then for
every zero $a$ of $p$, there exists a critical point $\zeta$
of $p$ such that $|a-\zeta|\le1$.
This long-standing conjecture was recently proved by Lech Mazur~\cite{Maz26},
and a streamlined exposition of the proof was subsequently
given by Tao \cite{Tao26}.
Two further strengthenings of Sendov's conjecture have also
been established: the quadratic Tang--Zhang inequality
\cite{Zha26b} and Borcea's $2$-variance conjecture
\cite{BF26,Zha26c}. 
In \cite{Zha26c}, the author also identified a common proof
strategy underlying these strengthenings and suggested
applying it to a broader class of strengthenings of
Sendov's conjecture:

\medskip
\emph{Direct arguments in low degrees, analytic estimates
in high degrees, and rigorous computer-assisted verification
for the remaining intermediate degrees}.
\medskip

\noindent This strategy was subsequently applied to resolve the
Erd\H{o}s--Herzog--Piranian conjecture on the length
of polynomial lemniscates \cite{Zha26d}.
In the present paper, we apply the same strategy
to Schmeisser's conjecture.

In 1977, Schmeisser \cite{Sch77} (see also \cite[Conjecture~18]{Tao26}) proposed the following
strengthening of Sendov's conjecture.

\begin{conj}[Schmeisser]\label{conj:Schmeisser}
Let $p$ be a polynomial of degree $n\ge2$ with all its zeros
in the closed unit disk. Denote by $H(p)$ the convex hull
of the zeros of $p$. Then for every point $\xi\in H(p)$,
there exists a critical point $\zeta$ of $p$ such that
\[
|\xi-\zeta|\le1.
\]
Moreover, the bound $1$ is sharp: for
$p(z)=c(z^n-e^{i\theta})$, where
$c\in\mathbb C\setminus\{0\}$ and $\theta\in\mathbb R$,
every zero of $p$ is at distance one from its only critical point.
\end{conj}

Taking $\xi$ to be any zero of $p$,
Conjecture~\ref{conj:Schmeisser} reduces to Sendov's conjecture. The known results concerning Conjecture~\ref{conj:Schmeisser}
are summarized in Table~\ref{tab:schmeisser}.

\begin{table}[htbp]
\centering
\small
\renewcommand{\arraystretch}{1.2}
\begin{tabular}{|c|p{7.0cm}|c|}
\hline
Upper bound & Hypotheses & Reference \\
\hline
$1$
& $n\le 4$
& Schmeisser \cite{Sch69} \\

$1$
& $n=5$
& Gács \cite{Gac71} \\

$1$
& Every vertex of $H(p)$ lies on the unit circle
& Schmeisser \cite{Sch77} \\

$1$
& $H(p)$ is a triangle
& Schmeisser \cite{Sch77} \\

$1$
& $p$ is a trinomial
& Schmeisser \cite{Sch71,Sch77} \\
\hline

$\displaystyle
\frac{2n^{1/(n-1)}}{n^{2/(n-1)}+1}$
& $\xi$ is the centroid of the zeros of $p$
& Paw{\l}owski \cite{Paw98} \\

$\displaystyle \sqrt{\frac{n-2}{n-1}}$
& $\xi$ is the centroid of the zeros of $p$
& Zhang \cite{Zha26a} \\

$\displaystyle \sqrt{\frac{n-2}{n-1}}\sqrt{1-|\xi|^2}$
& $\xi$ is the centroid of the zeros of $p$
& Zhang \cite{Zha27} \\

$\displaystyle \sqrt{1-|\xi|^2}$
& $\xi\in H(p')$
& Zhang \cite{Zha27} \\

$\displaystyle \sqrt{1-|\xi|^2}$
& Every extreme zero of $p$ is multiple
& Zhang \cite{Zha27} \\
\hline
\end{tabular}
\caption{Known cases and refinements of Conjecture~\ref{conj:Schmeisser}. Here $H(p)$ denotes the convex hull of the
zeros of $p$. The upper bounds are for the distance
from $\xi$ to the nearest critical point of $p$.
Unless otherwise specified, $\xi$ is an arbitrary
point of $H(p)$.}
\label{tab:schmeisser}
\end{table}

\begin{rem}\label{rem:stronger-disk}
Goodman, Rahman, and Ratti~\cite{GRR69} and, independently,
Schmeisser~\cite{Sch69} proposed the following conjecture.
If $p$ is a polynomial of degree $n\ge2$ with all its zeros
in the closed unit disk, then for each zero $a$ of $p$
there is a critical point $\zeta$ of $p$ such that
\[
\left|\zeta-\frac a2\right|\le1-\frac{|a|}{2}.
\]
This disk condition implies Conjecture~\ref{conj:Schmeisser}.
Indeed, for any $z\in H(p)$,
Lemma~\ref{lem:diameter-cover} gives a zero $a$ of $p$
with $|z-a/2|\le|a|/2$. The corresponding critical point
$\zeta$ then satisfies
\[
|z-\zeta|
\le \left|z-\frac a2\right|
   +\left|\zeta-\frac a2\right|
\le \frac{|a|}{2}+1-\frac{|a|}{2}
=1.
\]

The case $|a|=1$ holds in every degree; see Meir and
Sharma~\cite[Theorem~1]{Mei69}.
After rotating so that $a=1$, the conclusion is that $p'$
has a zero $\zeta$ with $|\zeta-\tfrac12|\le\tfrac12$.
Miller~\cite{Mil90} constructed counterexamples to this
conjecture in degrees $6$, $8$, $10$, and $12$.
Nevertheless, the geometric argument above suggests the
approach we take here: for each zero $a$, we cover the disk
with diameter $[0,a]$ by the unit disks centred at the
critical points of $p$, allowing the critical point to
depend on the point being covered.
\end{rem}

\begin{rem}
The author~\cite[Conjecture~2.4]{Zha26a} proposed a common
strengthening of Schmeisser's conjecture and Borcea's
variance conjectures, replacing the unit radius in
Conjecture~\ref{conj:Schmeisser} by the $r$-variance
\[
\sigma_r(p):=
\min_{c\in\mathbb C}
\left(\frac1n\sum_{j=1}^{n}|z_j-c|^r\right)^{1/r},
\qquad 1\le r<\infty,
\]
where $z_1,\ldots,z_n$ are the zeros of $p$, counted with
multiplicity. The conjectured bound is
\[
\min_{p'(\zeta)=0}|\xi-\zeta|\le\sigma_r(p)
\qquad\text{for every }\xi\in H(p).
\]
Taking $\xi$ to be a zero of $p$ gives Borcea's
$r$-variance conjecture. When all zeros lie in the closed
unit disk, the inequality $\sigma_r(p)\le1$ shows that
this bound also implies Schmeisser's conjecture.

This strengthening fails for every finite $r\ge1$.
To see this, fix $r$ and choose an integer $n>3^r$.
Let $p(z)=z^{n-1}(z-1)$. Then $H(p)=[0,1]$ and
$p'(z)=z^{n-2}(nz-(n-1))$, so the critical points are
$0$ and $(n-1)/n$.
At their midpoint $\xi=(n-1)/(2n)\in H(p)$, we have
\[
\min_{p'(\zeta)=0}|\xi-\zeta|
=\frac{n-1}{2n}
>\frac13
>n^{-1/r}
\ge\sigma_r(p),
\]
where the last inequality follows by taking $c=0$
in the definition of $\sigma_r(p)$.
\end{rem}

Let $p(z)=\prod_{j=1}^{n}(z-z_j)$, where $n\ge2$.
The incomplete polynomials associated with $p$ are
\[
g_j(z):=\prod_{\substack{1\le k\le n\\k\ne j}}(z-z_k),
\qquad 1\le j\le n.
\]
For $\lambda=(\lambda_1,\ldots,\lambda_n)\in[0,\infty)^n$
with $\sum_{j=1}^{n}\lambda_j=1$, write
$q_\lambda(z):=\sum_{j=1}^{n}\lambda_jg_j(z)$.
 The author \cite[Theorem~5.7]{Zha24} proved that
\[
H(p)=
\bigcup_{\substack{\lambda\in[0,\infty)^n\\
                   \sum_{j=1}^{n}\lambda_j=1}}
\{\xi\in\mathbb C:q_\lambda(\xi)=0\}.
\]
Thus, Conjecture~\ref{conj:Schmeisser} has the following consequence.

\begin{cor}\label{cor:incomplete-polynomials}
Let $p$ be a polynomial of degree $n\ge2$ with all its zeros
in the closed unit disk, and let $q_\lambda$ be any convex
combination of the incomplete polynomials associated with $p$.
Then for every zero $\xi$ of $q_\lambda$, there exists
a critical point $\zeta$ of $p$ such that
\[
|\xi-\zeta|\le1.
\]
\end{cor}

In this paper, we prove Conjecture~\ref{conj:Schmeisser} for every degree $n\ge2$. 
While polishing this manuscript, we became aware of two
recent preprints by Tang: 
one claiming to prove
Conjecture~\ref{conj:Schmeisser} for all sufficiently
large degrees \cite{Tan26a}, and the other claiming
to prove it for every $n\ge6251$ \cite{Tan26b}.
 The latter
manuscript runs to 101 pages!

Our proof rests on a relatively simple idea and differs
substantially from Tang's approach and from the methods
used in earlier proofs of Sendov's conjecture and its
stronger variants.
 We first cover the convex hull of the zeros by
 diameter disks and reduce the conjecture to the
 positive semidefiniteness of a $2\times2$ matrix.
 Using reflected polynomials and bounds for centred
 products, we then derive two necessary integral
 inequalities whose simultaneous validity would
 contradict the matrix condition. We exclude this
 possibility by analytic estimates in large degrees
 and rigorous computer-assisted verification in
 the remaining finite range.

\subsection{Notation}
We use the reciprocal moments from
Zhang~\cite[Section~2]{Zha26c}.
Let $p$ be a polynomial of degree $n\ge3$ with all its zeros in
$\D=\{z\in\mathbb C:|z|\le1\}$.
At a nonzero simple zero, after a rotation and multiplication
by a nonzero constant, we may write
\[
p(z)=(z-a)\prod_{j=1}^{n-1}(z-z_j),
\qquad 0<a\le1.
\]
Let $\zeta_1,\ldots,\zeta_{n-1}$ be the critical points of $p$,
counted with multiplicity. Since $p'(a)\ne0$, the following
quantities are well defined:
\begin{itemize}
\item
$q_j:=\dfrac{1}{a-\zeta_j}$ for $1\le j\le n-1$.

\item
$\mu:=\dfrac{1}{n-1}\sum_{j=1}^{n-1}q_j=x+iy$,
where $x,y\in\mathbb R$.

\item
$s:=\dfrac{1}{n-1}\sum_{j=1}^{n-1}|q_j|^2$.
\end{itemize}

We also set $\Delta_s:=(1-a^2)s+2ax-1$.
The Gauss--Lucas theorem~\cite[p.~18]{BE95} gives $|\zeta_j|\le1$ for every $j$.
Since $1-aq_j=-\zeta_jq_j$, we have
\begin{equation}\label{eq:delta-positive}
\Delta_s
=\frac{1}{n-1}\sum_{j=1}^{n-1}
(1-|\zeta_j|^2)|q_j|^2
\ge0.
\end{equation}

\subsection{High-level description of the proof}\label{sec:ideas}

Schmeisser~\cite[Section~3]{Sch77} observed that the stronger disk
condition recalled in Remark~\ref{rem:stronger-disk} implies the
convex-hull conclusion. The geometry of this implication suggests the
covering problem with which we begin.

\medskip
\noindent
\(\bullet\)\enspace\textbf{Step I. A weighted covering criterion.}
For a zero $a$ of $p$, let
$D_a=\{z\in\mathbb C:|z-a/2|\le|a|/2\}$.\footnote{The disk $D_a$ has
diameter $[0,a]$, and its defining inequality is equivalent to
$\Rea(z\overline a)\ge|z|^2$. For $a\ne0$ and $z\notin\{0,a\}$,
this says that the angle subtended by $[0,a]$ at $z$ is at least a
right angle.}
This disk has the same centre as the disk in
Remark~\ref{rem:stronger-disk}, but its radius is $|a|/2$.
 Thus, if a critical point $\zeta$ satisfies
$|\zeta-a/2|\le1-|a|/2$, then
$|z-\zeta|\le|z-a/2|+|\zeta-a/2|\le1$ for every $z\in D_a$.
Lemma~\ref{lem:diameter-cover} shows that the disks $D_a$, as $a$ ranges
over the zeros, cover $H(p)$. This explains why the stronger disk
condition implies the conjecture.

The counterexamples of Miller~\cite{Mil90} show that this stronger condition
fails in general. We therefore seek to cover each $D_a$ by the unit
disks centred at all the critical points. The critical point used to
cover $z$ may now depend on $z$.
At a nonzero simple zero, rotated to $a\in(0,1]$, 
 set
$\omega_j=|q_j|^2/((n-1)s)$.\footnote{The weights are proportional to
the inverse squared distances from $a$ to the critical points. Their
weighted centroid is $\sum_j\omega_j\zeta_j=a-\overline\mu/s$, and the
variance about it is $(s-|\mu|^2)/s^2$. The mean squared distance in
\eqref{eq:weighted} is therefore the squared distance from this
centroid plus its variance.}
Since these weights are positive and sum to one, it suffices to bound
the weighted mean of $|z-\zeta_j|^2$ by one on $D_a$.
Identity~\eqref{eq:weighted} expresses this mean in terms of $a$, $\mu$,
and $s$. Maximizing over $D_a$ gives the matrix condition
\[
\begin{pmatrix}s-1&ay\\ ay&\Delta_s\end{pmatrix}\ge0
\]
of Lemma~\ref{lem:geometry}.
Multiple zeros and the zero at the origin are dealt with directly.
At $a=1$, the reciprocal-sum argument of Meir and
Sharma~\cite[Section~2]{Mei69} gives the matrix condition; see
Lemma~\ref{lem:boundary-zero}. The remaining case is $0<a<1$.

\medskip
\noindent
\(\bullet\)\enspace\textbf{Step II. Choosing the comparison point.}
Suppose that \eqref{eq:matrix} fails. We keep $a$ and $\mu$ fixed and
replace $s$ by the smallest bound $\widehat s\ge\max\{1,s\}$ satisfying
$\Delta_{\widehat s}(\widehat s-1)=a^2y^2$, where
$\Delta_{\widehat s}=(1-a^2)\widehat s+2ax-1$.\footnote{With $s$
replaced by $\widehat s$, the matrix in \eqref{eq:matrix} is positive
semidefinite of rank one. Thus the quadratic form lies on the boundary
of the positive semidefinite cone. Only the bound for the second
moment changes; the critical points remain fixed.}
This equality leads us to define
\[
\gamma=a-\frac{\overline\mu}{\widehat s},\qquad
u=\frac{|\mu|^2}{\widehat s},\qquad
h=(1-a^2)(1-u+\Delta_{\widehat s}).
\]
The construction gives $|\gamma|<1$, $0<u\le1$, and
$0<h<1$.\footnote{Since
$\gamma=(1-s/\widehat s)a+(s/\widehat s)\sum_j\omega_j\zeta_j$,
the point $\gamma$ lies on the segment from $a$ to the weighted
centroid. Also, $u=\widehat s|a-\gamma|^2$ and
$h=(1-a^2)\widehat s(1-|\gamma|^2)$, so $u$ measures a scaled squared
displacement, while $h$ involves the deficits of the squared moduli
of $a$ and $\gamma$ from one.}
The choice of $\gamma$ is dictated by the factors
$1+(\gamma-a)tq_j$ in the normalized derivative
$p'(a+t(\gamma-a))/p'(a)$: their mean is the real number $1-ut$, and
their mean squared deviation is at most $u(1-u)t^2$.
We can therefore compare their product with $(1-ut)^{n-1}$.
The reflected second moment has the form \eqref{eq:E},
with the additional term $-ht(1-t)$.
The determinant equality gives $\Delta_{\widehat s}\le u$,
which is needed to compare the two endpoint integrals.

\medskip
\noindent
\(\bullet\)\enspace\textbf{Step III. Comparing a polynomial with its reflection.}
Put $g(z)=p(z)/(z-a)$ and
$g^*(z)=z^{n-1}\overline{g(1/\overline z)}$, extended polynomially
at zero.\footnote{Reflection sends each nonzero zero $z_j$ of $g$ to
$1/\overline{z_j}$. The polynomials $g$ and $g^*$ have the same modulus
on the unit circle. A zero of $g$ at the origin lowers the degree of
$g^*$, as inversion sends the origin to infinity.}
Since the zeros of $g$ lie in the closed unit disk,
$|g(z)|\le|g^*(z)|$ for $|z|\le1$.
This is the usual modulus bound for a finite Blaschke product: after
cancelling common factors on the unit circle, $g/g^*$ is such a
product, possibly a unimodular constant. For this interpretation, see
Garcia, Mashreghi, and Ross~\cite[Section~2]{GMR17}.
We prove the inequality directly by comparing the individual factors.

The fundamental theorem of calculus writes $g$ as the integral of a
product in \eqref{eq:integral-representation}; this representation
also appears in Mazur~\cite[Section~3]{Maz26}.
We apply the reflected inequality at $a$ and at the point $\gamma$
chosen in Step~II. At $a$, it yields the first inequality in
Lemma~\ref{lem:reduction}, which will give a lower bound for $h$.
At $\gamma$, we estimate the centred product from below and the
reflected product from above. The arithmetic--geometric mean
inequality and a convexity argument give the second inequality,
whose upper bound involves the positive part of
$1-2ut+ut^2-ht(1-t)$.
The subtraction of $ht(1-t)$ is essential: it must outweigh the error
in replacing the centred product by $(1-ut)^{n-1}$.

\medskip
\noindent
\(\bullet\)\enspace\textbf{Step IV. Estimating the centred product.}
We use two estimates for this error. The first rests on the
second-moment inequality conjectured by Schoenberg~\cite{Sco86} and
proved by Pereira~\cite[Section~3]{Per03} and
Malamud~\cite[Theorem~4.15]{Mal05}. Its centred form is stated in
Lemma~\ref{lem:schoenberg}.
Repeated differentiation, with this inequality applied at each stage,
controls the mean square of the zeros at each lower degree.
The arithmetic--geometric mean inequality then bounds the constant
term, giving the first coefficient bound in
\eqref{eq:coefficient-majorant}. We also obtain bounds from Cauchy's
coefficient estimate and Newton's identities, and use the smallest
of the three.

The second estimate follows by moving each factor along the segment
from its common mean to its actual value. At the initial
configuration, the derivative of the product vanishes because the
deviations sum to zero. Identity~\eqref{eq:centred-derivative}
expresses the remaining derivative in terms of the squares of these
deviations. Integrating gives Lemma~\ref{lem:interpolation} and then
\eqref{eq:error-integral}.
The coefficient bounds are particularly useful in small degrees and
near $u=1$, while the integral estimate is suited to the large-degree
argument. Both preserve the factor $1-u$ coming from the variance.

\medskip
\noindent
\(\bullet\)\enspace\textbf{Step V. Excluding the necessary inequalities.}
Lemma~\ref{lem:polar} gives
$u\ge((n+1)/2)^{2/(n-1)}-1$ and $(n-1)hu>6(1-u)$.
For fixed $u$, the right-hand side of \eqref{eq:second} decreases
with $h$. Replacing $h$ by a lower bound therefore leaves, for each
degree, a comparison in the single variable $u$.
Integration by parts gives the recurrence of
Lemma~\ref{lem:recurrence} for the quadratic integral. We combine its
upper bounds with the product estimates from Step~IV to obtain the
strict reverse of \eqref{eq:quadratic-comparison}.
As $u$ approaches one, both the gain and the error tend to zero.
The normalization in \eqref{eq:acceptance} removes their common
factor $1-u$, which is why retaining that factor in Step~IV matters.

The cases $n=2,3$ are direct. For $n\ge6001$, the analytic estimates
in Lemma~\ref{lem:large} prove the required strict inequality
uniformly in $u$. For $4\le n\le6000$, Lemma~\ref{lem:finite} uses
interval certificates for the admissible range of $u$.
The verifier checks both the inequalities, using outward-rounded
integer arithmetic, and the coverage by closed intervals, including
their shared endpoints.

\subsection{Organization of the paper}

Section~\ref{sec:geometry} proves the covering criterion in Step~I.
Section~\ref{sec:reduction} carries out Steps~II and~III and derives the
scalar parameter bounds. Section~\ref{sec:error} proves the two product
estimates in Step~IV, and Section~\ref{sec:moments} develops the quadratic
integral estimates used in Step~V. Section~\ref{sec:exclusion} excludes
the necessary inequalities and completes the proof of
Conjecture~\ref{conj:Schmeisser}. Appendix~\ref{app:verification} gives the
endpoint estimates, interval bounds, arithmetic rules, and reproduction
instructions for the finite verification.

\subsection{Formalization in Lean 4}\label{sec:lean}

We have formalized the proof of
Conjecture~\ref{conj:Schmeisser} for every degree $n\ge2$,
together with the optimality of the constant $1$ in each
degree, in Lean~4 using mathlib.
The source files, certificates, and build instructions
are available in the accompanying repository~\cite{Zha26e}.

The main result and the optimality statement appear as
\texttt{Schmeisser.main} and
\texttt{Schmeisser.universal\_bound\_iff}, respectively,
in \texttt{lean/Schmeisser/Main.lean}.
The formal proof includes verification of the finite
interval certificates by Lean's kernel.
It contains no \texttt{sorry} placeholders and uses no
axioms beyond Lean's standard logical foundations.

\subsection{Acknowledgments and AI tools disclosure} The author thanks Alexandre Eremenko for his valuable
suggestions, Minghua Lin for suggesting the use of matrix
methods to study the geometry of polynomials, and
Terence Tao for very helpful discussions in the comments
on his
\href{https://terrytao.wordpress.com/2026/08/12/a-digestion-of-the-proof-of-sendovs-conjecture/}{blog post}.

This work was supported by the China Scholarship Council,
by the Young Elite Scientists Sponsorship Program for PhD Students
of the China Association for Science and Technology,
and by the Fundamental Research Funds for the Central Universities
at Xi'an Jiaotong University (Grant No.~xzy022024045).

ChatGPT 6.1 was used for English-language editing, proofreading, grammatical corrections, the preparation of illustrative figures, and exploratory discussions of possible approaches to selected parts of this paper. Throughout this process, the author provided mathematical guidance, critically evaluated and corrected errors, and adjusted the research directions as necessary. The mathematical arguments and proofs in the original manuscript were independently developed, checked, and written by the author.

The Lean~4 formalization in~\cite{Zha26e} was generated using
OpenAI's Codex and checked by the Lean kernel.

\section{A weighted covering criterion}\label{sec:geometry}

Recall for $a\in \mathbb C$,
\[
D_a=\left\{z\in\mathbb C:\left|z-\frac a2\right|\le\frac{|a|}{2}\right\}.
\]

\begin{lem}[Covering by diameter disks]\label{lem:diameter-cover}
Let $p$ be a polynomial of degree $n\ge1$ with zeros
$z_1,\ldots,z_n$. Then
\[
H(p)\subseteq\bigcup_{j=1}^n D_{z_j}.
\]
\end{lem}

\begin{proof}
Let $z\in H(p)$. Then
\[
z=\sum_{j=1}^n t_jz_j,\qquad
 t_j\ge0,\qquad\sum_{j=1}^n t_j=1.
\]
Consequently,
\[
\sum_{j=1}^n t_j\Rea(z\overline{z_j})
=\Rea\left(z\overline{\sum_{j=1}^n t_jz_j}\right)
=|z|^2.
\]
Thus $\Rea(z\overline{z_j})\ge|z|^2$ for some $j$. This is equivalent to
$|z-z_j/2|\le|z_j|/2$, and hence $z\in D_{z_j}$.
\end{proof}

\begin{lem}[A weighted covering criterion]\label{lem:geometry}
Let $a\in(0,1]$ be a simple zero of $p$. If
\begin{equation}\label{eq:matrix}
\mathsf M(s):=\begin{pmatrix}s-1&ay\\ay&\Delta_s\end{pmatrix}\ge0,
\end{equation}
then every point of $D_a$ is within distance one of a critical point of $p$.
\end{lem}

\begin{proof}
Set
\[
\omega_j=\frac{|q_j|^2}{(n-1)s},\qquad 1\le j\le n-1.
\]
These weights are positive and sum to one. Since
$\zeta_j=a-1/q_j$, a direct calculation gives
\begin{equation}\label{eq:weighted}
\sum_{j=1}^{n-1}\omega_j|z-\zeta_j|^2
=|z-a|^2+\frac{2\Rea((z-a)\mu)+1}{s}.
\end{equation}
The expression on the right is convex in $z$, so its maximum over $D_a$
is attained on the boundary. Writing $z=a(1+e^{it})/2$, we obtain
\begin{equation}\label{eq:weighted-max}
\max_{z\in D_a}\sum_{j=1}^{n-1}\omega_j|z-\zeta_j|^2
=\frac{a^2}{2}+\frac{1-ax}{s}
 +\frac{a}{s}\left|\mu-\frac{as}{2}\right|.
\end{equation}
By the definition of $\Delta_s$, the right-hand side of
\eqref{eq:weighted-max} is at most one if and only if
\begin{equation}\label{eq:weighted-condition}
a\left|\mu-\frac{as}{2}\right|
\le\frac{s-1+\Delta_s}{2}.
\end{equation}
The assumption \eqref{eq:matrix} implies $s-1\ge0$ and
$\Delta_s\ge0$. Furthermore,
\begin{equation}\label{eq:weighted-determinant}
\left(\frac{s-1+\Delta_s}{2}\right)^2
-a^2\left|\mu-\frac{as}{2}\right|^2
=(s-1)\Delta_s-a^2y^2\ge0.
\end{equation}
Since the right-hand side of \eqref{eq:weighted-condition} is
nonnegative, \eqref{eq:weighted-determinant} proves
\eqref{eq:weighted-condition}. Thus the weighted mean in
\eqref{eq:weighted} is at most one for every $z\in D_a$.
As the weights sum to one, at least one critical point $\zeta_j$
satisfies $|z-\zeta_j|\le1$.
\end{proof}

\begin{figure}[tbp]
\centering
\includegraphics[width=\textwidth]{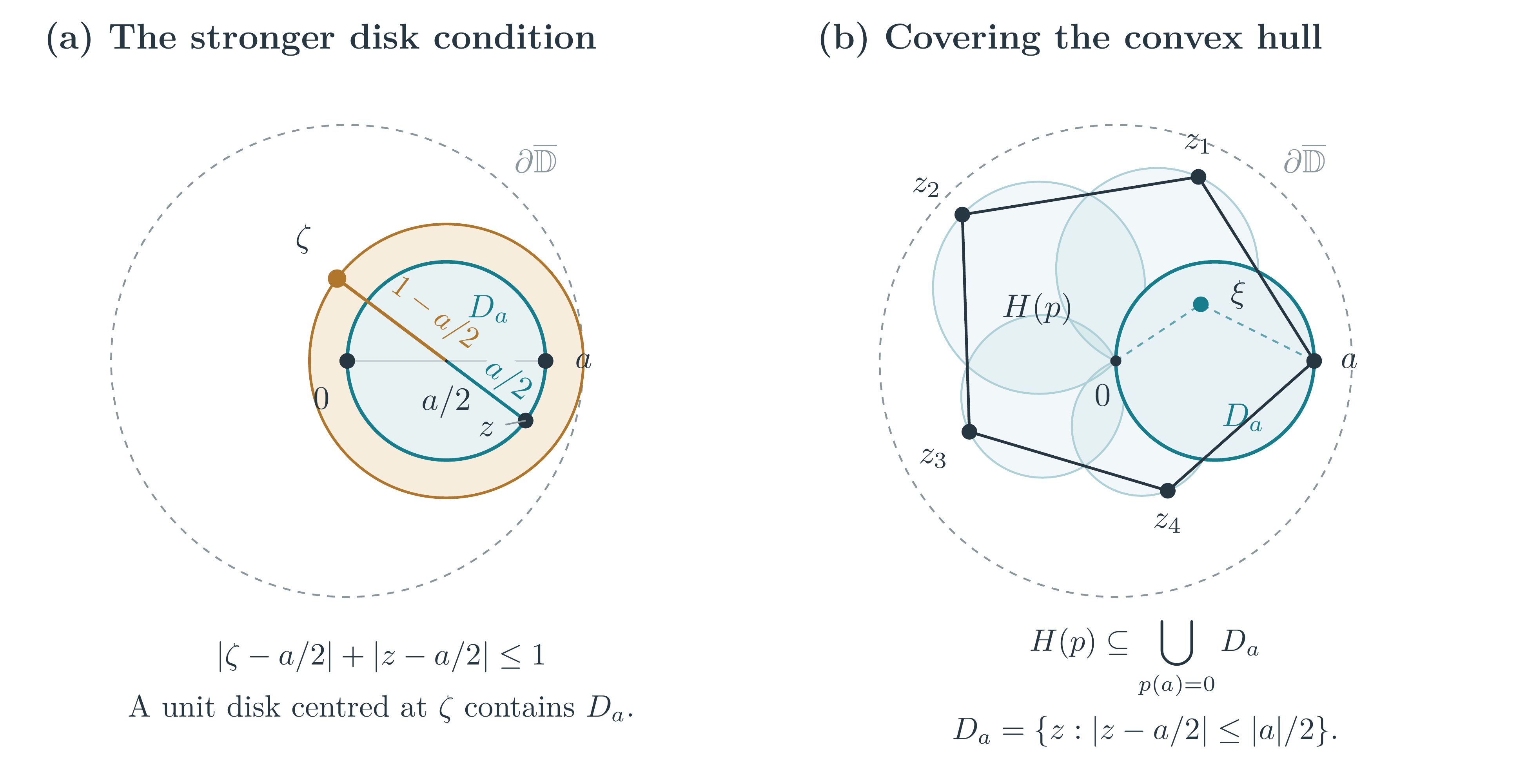}
\caption{The geometric reduction.
(a) For $0<a<1$, the condition
$|\zeta-a/2|\le 1-a/2$ ensures that the unit disk centred
at $\zeta$ contains $D_a$; see Remark~\ref{rem:stronger-disk}.
(b) The diameter disks associated with the zeros of $p$
cover $H(p)$, as in Lemma~\ref{lem:diameter-cover}.
The darker disk is $D_a$, and $\xi\in H(p)\cap D_a$.}
\label{fig:diameter}
\end{figure}

We now explain the reduction to \eqref{eq:matrix}. Given any
$z\in H(p)$, Lemma~\ref{lem:diameter-cover} supplies a zero $z_i$
with $z\in D_{z_i}$; see Figure~\ref{fig:diameter}.
If $z_i$ is a multiple zero, it is a critical point, and
\[
|z-z_i|\le\left|z-\frac{z_i}{2}\right|+\frac{|z_i|}{2}
\le|z_i|\le1.
\]
If $z_i=0$, then $D_{z_i}=\{0\}$, so $z=0$; the Gauss--Lucas theorem
implies that every critical point lies in the unit disk.
It remains to consider a nonzero simple zero $z_i$. After a rotation,
we may assume $z_i=a\in(0,1]$. By Lemma~\ref{lem:geometry}, it
suffices to establish \eqref{eq:matrix} at every such zero.

The matrix condition \eqref{eq:matrix} is automatic at a simple zero
on the unit circle.

\begin{lem}[A zero on the unit circle]\label{lem:boundary-zero}
If $a=1$ is a simple zero of $p$, then \eqref{eq:matrix} holds.
\end{lem}

\begin{proof}
Write
\[
p(z)=(z-1)\prod_{j=1}^{n-1}(z-z_j),
\]
where $|z_j|\le1$ and $z_j\ne1$ for $1\le j\le n-1$.
By the Meir--Sharma identity~\cite{Mei69},
\begin{equation}\label{eq:meir-sharma}
\frac{p''(1)}{p'(1)}
=\sum_{j=1}^{n-1}\frac{1}{1-\zeta_j}
=2\sum_{j=1}^{n-1}\frac{1}{1-z_j}.
\end{equation}
Taking real parts in \eqref{eq:meir-sharma} and recalling the definition
of $\mu=x+iy$ gives
\begin{equation}\label{eq:boundary-x}
x=\frac{2}{n-1}\sum_{j=1}^{n-1}
\Rea\frac{1}{1-z_j}\ge1.
\end{equation}
Indeed,
\[
\Rea\frac{1}{1-z_j}-\frac12
=\frac{1-|z_j|^2}{2|1-z_j|^2}\ge0
\qquad(1\le j\le n-1).
\]
Since $a=1$, the definition of $\Delta_s$ and
\eqref{eq:boundary-x} give $\Delta_s=2x-1\ge1$.
By the Cauchy--Schwarz inequality,
\[
s=\frac{1}{n-1}\sum_{j=1}^{n-1}|q_j|^2
\ge\left|\frac{1}{n-1}\sum_{j=1}^{n-1}q_j\right|^2
=x^2+y^2.
\]
As $x\ge1$, it follows that $s-1\ge y^2\ge0$. Consequently,
\[
(s-1)\Delta_s-y^2\ge(s-1)-y^2\ge0.
\]
Both diagonal entries and the determinant of the matrix in
\eqref{eq:matrix} are therefore nonnegative, proving the assertion.
\end{proof}

\section{Two necessary integral inequalities}\label{sec:reduction}

For $n\ge3$, $0<u\le1$, and $0\le h\le1$, write
$\eps=1-u$\footnote{In the construction of
Lemma~\ref{lem:reduction}, $\eps$ is the proportion of the enlarged
second moment left after subtracting the squared modulus of the mean.
The identity $\mathcal B_u(t)=(1-ut)^2+u\eps t^2$ expresses the quadratic
below as a squared Euclidean norm. For centred perturbations
$\eta_1,\ldots,\eta_{n-1}$ with mean square at most $u\eps$, it bounds
the mean squared modulus of the factors $1-ut+t\eta_j$.} and
\[
\mathcal B_u(t)=1-2ut+ut^2,
\qquad 0\le t\le1.
\]
Define
\begin{equation}\label{eq:C}
C_{n-1}(u):=
n\int_0^1(1-ut)^{n-1}\dd t
=\frac{1-\eps^n}{u}
=\sum_{j=0}^{n-1}\eps^j
\end{equation}
and
\[
P_{n-1}(u,h):=
\int_0^1[1-h+h(1+u)t]^{(n-1)/2}\dd t.
\]
We also introduce the centred-product error.\footnote{The constraints
below describe a configuration of $n-1$ points with centroid zero and
mean squared distance from the origin at most $u\eps$. In the proof of
Lemma~\ref{lem:reduction}, these points are
$\eta_j=-\overline\mu(q_j-\mu)/\widehat s$: a translation followed by a
rotation and dilation of the reciprocal critical points.}
\begin{equation}\label{eq:error-definition}
\Err_{n-1}(u):=
n\max_{\substack{
\eta_1,\ldots,\eta_{n-1}\in\mathbb C,\ \sum_j\eta_j=0\\
\sum_j|\eta_j|^2\le(n-1)u\eps}}
\int_0^1
\left|
\prod_{j=1}^{n-1}(1-ut+t\eta_j)
-(1-ut)^{n-1}
\right|\dd t.
\end{equation}
The maximum is attained by compactness and continuity.
In particular, $\Err_{n-1}(1)=0$, since the constraint
then forces every $\eta_j$ to vanish.
For a real number $r$, we write $r_+=\max\{r,0\}$.

\begin{lem}[Necessary integral inequalities]
\label{lem:reduction}
If \eqref{eq:matrix} fails at a simple zero
$a\in(0,1)$, then there exist $0<u\le1$ and
$0<h<1$ such that
\begin{align}
1&\le P_{n-1}(u,h),\label{eq:polar}\\
C_{n-1}(u)-\Err_{n-1}(u)
&\le n\int_0^1
[\mathcal B_u(t)-ht(1-t)]_+^{(n-1)/2}\dd t.
\label{eq:second}
\end{align}
\end{lem}

\begin{proof}
\emph{Step 1. Enlarging the second moment.}
Fix a simple zero $a\in(0,1)$ at which
\eqref{eq:matrix} fails. Recall from
\eqref{eq:delta-positive} that $\Delta_s\ge0$.

For $\sigma\ge s$, define
$
\Delta_\sigma=(1-a^2)\sigma+2ax-1$.
Since $1-a^2>0$, we have
$\Delta_\sigma\ge\Delta_s\ge0$.
Consider the function
\[
\Phi(\sigma)
=(\sigma-1)\Delta_\sigma-a^2y^2,
\qquad \sigma\ge\max\{1,s\}.
\]

We first show that
$\Phi(\max\{1,s\})\le0$.
If $s\ge1$, both diagonal entries of the
matrix in \eqref{eq:matrix} are nonnegative.
Since this matrix is not positive semidefinite,
its determinant must be negative. Hence
$\Phi(s)<0$.
If $s<1$, then
$\Phi(1)=-a^2y^2\le0$.

Moreover,
\[
\Phi'(\sigma)
=(1-a^2)(\sigma-1)+\Delta_\sigma\ge0
\]
for $\sigma\ge\max\{1,s\}$, and
$\Phi(\sigma)\to\infty$ as $\sigma\to\infty$.
By continuity, there is a smallest
$\widehat s\ge\max\{1,s\}$ satisfying
\begin{equation}\label{eq:boundary}
\Delta_{\widehat s}(\widehat s-1)=a^2y^2.
\end{equation}

We claim that $\Delta_{\widehat s}>0$.
If $s<1$, then
\[
\Delta_{\widehat s}
\ge\Delta_1
=\Delta_s+(1-a^2)(1-s)>0.
\]
If $s\ge1$, the failure of \eqref{eq:matrix}
forces $y\ne0$, since otherwise the matrix
would be positive semidefinite.
In this case, \eqref{eq:boundary} implies
$\Delta_{\widehat s}>0$.

We now define
\[
u:=\frac{|\mu|^2}{\widehat s},
\qquad
h:=(1-a^2)(1-u+\Delta_{\widehat s}),
\qquad
\gamma:=a-\frac{\overline\mu}{\widehat s}.
\]
By the Cauchy--Schwarz inequality,
$|\mu|^2\le s\le\widehat s$, so $u\le1$.
To see that $u>0$, suppose that $\mu=0$.
Then $x=y=0$, and
$\Delta_s=(1-a^2)s-1\ge0$ implies $s>1$.
The matrix in \eqref{eq:matrix} would
therefore be positive semidefinite, a
contradiction.

We next establish the relation between
$u$ and $\Delta_{\widehat s}$.
By the definition of $\Delta_{\widehat s}$,
\[
2ax=1+\Delta_{\widehat s}
-(1-a^2)\widehat s.
\]
Combining this identity with
\eqref{eq:boundary} and
$u=(x^2+y^2)/\widehat s$, we obtain
\begin{equation}\label{eq:u-delta}
u-\Delta_{\widehat s}
=
\frac{
\bigl(1-\Delta_{\widehat s}
-(1-a^2)\widehat s\bigr)^2
}{4a^2\widehat s}
\ge0.
\end{equation}
Consequently,
\[
0<\Delta_{\widehat s}\le u\le1.
\]
It follows that
$0<1-u+\Delta_{\widehat s}\le1$.
Since $0<a<1$, the definition of $h$
now gives $0<h<1$.

Finally, direct calculation yields
\[
1-|\gamma|^2
=\frac{1-u+\Delta_{\widehat s}}{\widehat s}
>0.
\]
Thus $\gamma$ lies in the open unit disk.

\medskip
\noindent
\emph{Step 2. Comparison with the reflected polynomial.}
Let
$
g(z)={p(z)}/{(z-a)}.
$
Since $p(a)=0$, the fundamental theorem of
calculus gives
\[
g(\xi)=\int_0^1
p'(a+t(\xi-a))\dd t,
\]
where the identity at $\xi=a$ is understood
by continuity. Since the critical points
of $p$ are $\zeta_1,\ldots,\zeta_{n-1}$,
we have
\[
p'(z)=p'(a)\prod_{j=1}^{n-1}
\bigl(1+(z-a)q_j\bigr).
\]
Substitution into the preceding integral
gives
\begin{equation}\label{eq:integral-representation}
g(\xi)=p'(a)\int_0^1
\prod_{j=1}^{n-1}
\bigl(1+(\xi-a)tq_j\bigr)\dd t.
\end{equation}

Consider the reflected polynomial
\[
g^*(\xi)
=\xi^{n-1}\overline{g(1/\overline\xi)},
\]
with its polynomial extension at $\xi=0$.
Since $a$ is real, \eqref{eq:integral-representation}
gives
\begin{equation}\label{eq:reflected-integral}
g^*(\xi)
=\overline{p'(a)}
\int_0^1
\prod_{j=1}^{n-1}
\bigl[\xi+(1-a\xi)t\overline{q_j}\bigr]
\dd t.
\end{equation}

Every zero $z_j$ of $g$ lies in the closed unit disk.
Hence, for $|\xi|\le1$,
\[
|1-\overline{z_j}\xi|^2-|\xi-z_j|^2
=(1-|z_j|^2)(1-|\xi|^2)\ge0.
\]
Taking products over all zeros of $g$, we obtain
$|g(\xi)|\le|g^*(\xi)|$.
Substituting \eqref{eq:integral-representation}
and \eqref{eq:reflected-integral} into this inequality,
taking the complex conjugate of the integral
representing $g^*$, and cancelling the common
factor $|p'(a)|>0$, we obtain
\begin{equation}\label{eq:reflection}
\left|
\int_0^1
\prod_{j=1}^{n-1}
[1+(\xi-a)tq_j]\dd t
\right|
\le
\left|
\int_0^1
\prod_{j=1}^{n-1}
[\overline\xi+(1-a\overline\xi)tq_j]
\dd t
\right|,
\qquad |\xi|\le1.
\end{equation}

We first apply \eqref{eq:reflection}
with $\xi=\gamma$. By the definition of
$\gamma$ and $u$,
\[
1+(\gamma-a)tq_j
=1-ut-\frac{\overline\mu}{\widehat s}
(q_j-\mu)t.
\]
The coefficients
$-\overline\mu(q_j-\mu)/\widehat s$
have mean zero. Moreover, their mean
squared modulus is
\[
\frac{|\mu|^2}{\widehat s^2}
\frac1{n-1}\sum_{j=1}^{n-1}|q_j-\mu|^2=
\frac{|\mu|^2}{\widehat s^2}
(s-|\mu|^2)
\le u(1-u)=u\eps,
\]
where the inequality follows from
$s\le\widehat s$.
Thus these coefficients satisfy the
constraints in \eqref{eq:error-definition}.
The triangle inequality and \eqref{eq:C}
therefore show that the left-hand side
of \eqref{eq:reflection}, multiplied
by $n$, is at least
$C_{n-1}(u)-\Err_{n-1}(u)$.

Figure~\ref{fig:comparison-centering} illustrates the location
of $\gamma$ and the centred factors, with
$\eta_j=-\overline\mu(q_j-\mu)/\widehat s$.

\begin{figure}[tbp]
\centering
\includegraphics[width=\textwidth]{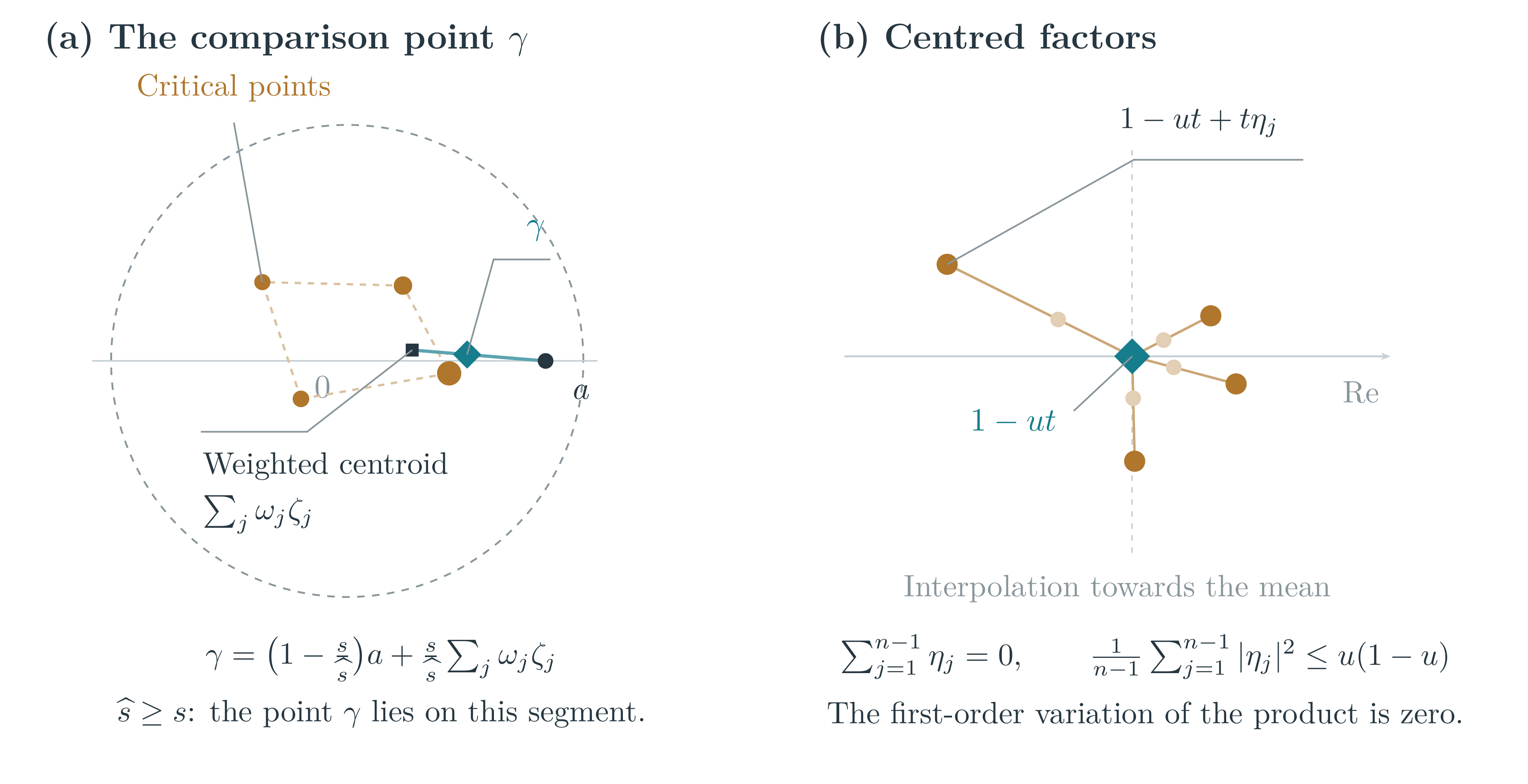}
\caption{The comparison point and centred interpolation.
(a) The point $\gamma$ lies on the segment joining $a$ to
$\sum_j\omega_j\zeta_j$.
The affine relation is illustrated with $\widehat s=1.7s$,
without imposing~\eqref{eq:boundary}.
(b) The factors $1-ut+rt\eta_j$, $0\le r\le1$, have
the same mean $1-ut$. Here $t=1$; pale and dark circles
correspond to $r=0.4$ and $r=1$, respectively.
Since $\sum_j\eta_j=0$, the derivative of their product
with respect to $r$ vanishes at $r=0$.}
\label{fig:comparison-centering}
\end{figure}

We next estimate the right-hand side of
\eqref{eq:reflection} with $\xi=\gamma$.
By the power-mean inequality and the definitions
of $\mu$ and $s$, for any $A,B\in\mathbb C$, we have
\begin{equation}\label{eq:product-rms}
\begin{aligned}
\left|\prod_{j=1}^{n-1}(A+Bq_j)\right|
&\le
\left(\frac{1}{n-1}\sum_{j=1}^{n-1}
|A+Bq_j|^2\right)^{(n-1)/2}\\
&=
\left(
|A+B\mu|^2+|B|^2(s-|\mu|^2)
\right)^{(n-1)/2}.
\end{aligned}
\end{equation}

Set
\[
\lambda:=1-|\gamma|^2
=\frac{1-u+\Delta_{\widehat s}}{\widehat s}.
\]
The preceding estimates show that
$0<\lambda\le1$.\footnote{Here $\lambda=1-|\gamma|^2$ measures the
deficit of $|\gamma|^2$ from one. The Euclidean distance from $\gamma$ to the unit circle is
$1-|\gamma|=\lambda/(1+|\gamma|)$.}
Applying \eqref{eq:product-rms} with
$A=\overline\gamma$ and
$B=(1-a\overline\gamma)t$, and then
using $s\le\widehat s$, gives the
nonnegative quadratic\footnote{Before $s$ is replaced by
$\widehat s$, the first line
is the mean squared modulus of the reflected factors
$\overline\gamma+(1-a\overline\gamma)tq_j$. The replacement enlarges
their variance while keeping their centroid fixed. Thus $E(t)$ bounds
their mean squared distance from the origin.}
\begin{align}
E(t)
&:=|\overline\gamma
+(1-a\overline\gamma)\mu t|^2
+|1-a\overline\gamma|^2
(\widehat s-|\mu|^2)t^2
\notag\\
&=(1-t)(1-\lambda)
+t(\eps+\Delta_{\widehat s}\lambda)
-(u+h)t(1-t).
\label{eq:E}
\end{align}

Combining \eqref{eq:C}, \eqref{eq:error-definition},
\eqref{eq:reflection}, \eqref{eq:product-rms},
and \eqref{eq:E}, we obtain
\begin{align}
C_{n-1}(u)-\Err_{n-1}(u)
&\le n\left|
\int_0^1\prod_{j=1}^{n-1}
[1+(\gamma-a)tq_j]\dd t
\right|\notag\\
&\le n\left|
\int_0^1\prod_{j=1}^{n-1}
[\overline\gamma+(1-a\overline\gamma)tq_j]
\dd t
\right|\notag\\
&\le n\int_0^1
\left(
|\overline\gamma+(1-a\overline\gamma)\mu t|^2
+|1-a\overline\gamma|^2
(s-|\mu|^2)t^2
\right)^{(n-1)/2}\dd t
\notag\\
&\le n\int_0^1 E(t)^{(n-1)/2}\dd t.
\label{eq:E-integral}
\end{align}

We next apply \eqref{eq:reflection} with $\xi=a$.
Since the integrand on the left-hand side is
identically one, the corresponding integral equals one.

To estimate the right-hand side, define
\begin{equation}\label{eq:H-definition}
\begin{aligned}
H(t)
&:=|a+(1-a^2)\mu t|^2
 +(1-a^2)^2(\widehat s-|\mu|^2)t^2\\
&=a^2+(1-a^2)(1+\Delta_{\widehat s})t
 -(1-a^2)^2\widehat s\,t(1-t).
\end{aligned}
\end{equation}
The second equality follows from
$\Delta_{\widehat s}
=(1-a^2)\widehat s+2ax-1$.
Since $\widehat s\ge s\ge|\mu|^2$,
the first expression in \eqref{eq:H-definition}
also shows that $H(t)\ge0$.\footnote{The same interpretation
applies at the point $a$:
$H(t)$ bounds the mean squared modulus of
$a+(1-a^2)tq_j$, with the actual second moment enlarged to $\widehat s$.}

Combining \eqref{eq:reflection} with
$\xi=a$, the triangle inequality, and
\eqref{eq:product-rms} with
$A=a$ and $B=(1-a^2)t$, we obtain
\begin{align}
1
&\le
\left|\int_0^1
\prod_{j=1}^{n-1}
[a+(1-a^2)tq_j]\dd t\right|
\notag\\
&\le
\int_0^1
\left|\prod_{j=1}^{n-1}
[a+(1-a^2)tq_j]\right|\dd t
\notag\\
&\le
\int_0^1
\left(
|a+(1-a^2)\mu t|^2
+(1-a^2)^2(s-|\mu|^2)t^2
\right)^{(n-1)/2}\dd t
\notag\\
&\le
\int_0^1H(t)^{(n-1)/2}\dd t.
\label{eq:H}
\end{align}
Here the last inequality follows from
$s\le\widehat s$ and \eqref{eq:H-definition}.

\medskip
\noindent
\emph{Step 3. Derivation of the integral inequalities.}

We first establish \eqref{eq:polar}.
Since $0\le t\le1$, the second expression
in \eqref{eq:H-definition} gives
\[
H(t)\le
a^2+(1-a^2)(1+\Delta_{\widehat s})t.
\]
On the other hand, the definition of $h$ yields
\[
1-h+h(1+u)t-
\bigl[a^2+(1-a^2)
(1+\Delta_{\widehat s})t\bigr]=(1-a^2)(u-\Delta_{\widehat s})(1-ut)
\ge0,
\]
where the last inequality follows from
\eqref{eq:u-delta} and $0<u\le1$.
Consequently,
\begin{equation}\label{eq:H-majorant}
0\le H(t)\le1-h+h(1+u)t,
\qquad 0\le t\le1.
\end{equation}
Combining \eqref{eq:H} and
\eqref{eq:H-majorant}, we obtain
\[
\begin{aligned}
1
&\le\int_0^1H(t)^{(n-1)/2}\dd t\\
&\le\int_0^1
[1-h+h(1+u)t]^{(n-1)/2}\dd t\\
&=P_{n-1}(u,h),
\end{aligned}
\]
which proves \eqref{eq:polar}.

We now turn to \eqref{eq:second}.
For $0\le\lambda\le1$, consider
\[
E_\lambda(t)
=(1-t)(1-\lambda)
+t(\eps+\Delta_{\widehat s}\lambda)
-(u+h)t(1-t).
\]
For each fixed $t$, this expression
is affine in $\lambda$. In particular,
\begin{equation}\label{eq:E-endpoints}
\begin{aligned}
E_0(t)
&=\mathcal B_u(t)-ht(1-t),\\
E_0(t)-E_1(1-t)
&=\eps t+(u-\Delta_{\widehat s})(1-t)
\ge0,
\end{aligned}
\end{equation}
where the last inequality follows from
\eqref{eq:u-delta}.

Since $r\mapsto r_+^{(n-1)/2}$ is
nondecreasing, \eqref{eq:E-endpoints}
and the change of variable $t\mapsto1-t$
give
\begin{equation}\label{eq:E-endpoint-integrals}
\int_0^1[E_1(t)]_+^{(n-1)/2}\dd t
=\int_0^1[E_1(1-t)]_+^{(n-1)/2}\dd t
\le\int_0^1[E_0(t)]_+^{(n-1)/2}\dd t.
\end{equation}

Moreover, since $(n-1)/2\ge1$, the function
$r\mapsto r_+^{(n-1)/2}$ is convex.
Using
\[
E_\lambda(t)
=(1-\lambda)E_0(t)+\lambda E_1(t),
\]
we obtain
\begin{align}
\int_0^1[E_\lambda(t)]_+^{(n-1)/2}\dd t
&\le
(1-\lambda)
\int_0^1[E_0(t)]_+^{(n-1)/2}\dd t
\notag+
\lambda
\int_0^1[E_1(t)]_+^{(n-1)/2}\dd t
\notag\\
&\le
\int_0^1[E_0(t)]_+^{(n-1)/2}\dd t,
\label{eq:E-convex-bound}
\end{align}
where the last inequality follows from
\eqref{eq:E-endpoint-integrals}.

For the value
$\lambda=1-|\gamma|^2\in(0,1]$
defined in Step~2, equation~\eqref{eq:E}
gives $E_\lambda(t)=E(t)\ge0$.
Therefore, combining \eqref{eq:E-integral},
\eqref{eq:E-convex-bound}, and
\eqref{eq:E-endpoints}, we conclude that
\[
\begin{aligned}
C_{n-1}(u)-\Err_{n-1}(u)
&\le n\int_0^1 E(t)^{(n-1)/2}\dd t\\
&\le n\int_0^1
[E_0(t)]_+^{(n-1)/2}\dd t\\
&=n\int_0^1
[\mathcal B_u(t)-ht(1-t)]_+^{(n-1)/2}\dd t.
\end{aligned}
\]
This is precisely \eqref{eq:second},
and the proof is complete.
\end{proof}

The inequality \eqref{eq:polar} imposes quantitative
lower bounds on both $u$ and $hu$.

\begin{lem}[Lower bounds for the scalar parameters]
\label{lem:polar}
Let $n\ge3$, $0<u\le1$, and $0<h<1$ satisfy
\eqref{eq:polar}. Then
\begin{equation}\label{eq:polar-bounds}
u\ge\tau_{n-1}:=
\left(\frac{n+1}{2}\right)^{2/(n-1)}-1,
\qquad
(n-1)hu>6(1-u).
\end{equation}
\end{lem}

\begin{proof}
We first establish the lower bound for $u$.
For fixed $u\in(0,1]$, the function
\[
h\longmapsto P_{n-1}(u,h)
=\int_0^1
[1-h+h(1+u)t]^{(n-1)/2}\dd t
\]
is convex on $[0,1]$, because $(n-1)/2\ge1$
and the integrand is a nonnegative affine
function of $h$ raised to this power.
Consequently,
\begin{equation}\label{eq:polar-chord}
\begin{aligned}
P_{n-1}(u,h)
&\le(1-h)P_{n-1}(u,0)
+hP_{n-1}(u,1)\\
&=1-h+\frac{2h}{n+1}(1+u)^{(n-1)/2}.
\end{aligned}
\end{equation}
Combining \eqref{eq:polar} and
\eqref{eq:polar-chord}, and using $h>0$,
we obtain
\[
(1+u)^{(n-1)/2}\ge\frac{n+1}{2}.
\]
Taking the power $2/(n-1)$ proves the
first inequality in \eqref{eq:polar-bounds}.

To establish the lower bound for $hu$,
set
\[
\ell=\frac{(n-1)h}{2},
\qquad
z=\frac{\ell(1+u)}{2}.
\]
Both $\ell$ and $z$ are positive.
Since $0<h<1$, we have
$1-h+h(1+u)t>0$ for $0\le t\le1$.
Applying $\log(1+r)<r$ for $r>-1$
and $r\ne0$ gives
\[
[1-h+h(1+u)t]^{(n-1)/2}
\le\exp\bigl(\ell((1+u)t-1)\bigr),
\]
with strict inequality except at
$t=1/(1+u)$. Integrating over $[0,1]$
and using \eqref{eq:polar}, we obtain
\begin{equation}\label{eq:hyperbolic}
\begin{aligned}
1
&\le P_{n-1}(u,h)\\
&<\int_0^1
\exp\bigl(\ell((1+u)t-1)\bigr)\dd t\\
&=\exp\left(-\frac{\ell(1-u)}{2}\right)
\frac{\sinh z}{z}.
\end{aligned}
\end{equation}

We use the following inequality established in
\cite[(27)]{Tao26}:
\begin{equation}\label{eq:hyperbolic-majorant}
\log\frac{\sinh z}{z}
\le\sqrt{z^2+9}-3,
\qquad z>0.
\end{equation}
Combining \eqref{eq:hyperbolic} and
\eqref{eq:hyperbolic-majorant}, we obtain
\[
\frac{\ell(1-u)}{2}
<\sqrt{z^2+9}-3.
\]
Since $z=\ell(1+u)/2$, squaring both sides
after adding $3$ gives
\[
3\ell(1-u)
<\frac{\ell^2}{4}
\bigl((1+u)^2-(1-u)^2\bigr)
=\ell^2u.
\]
Recalling that $\ell=(n-1)h/2>0$, we conclude that
\[
(n-1)hu>6(1-u),
\]
which proves the second inequality in
\eqref{eq:polar-bounds}.
\end{proof}

The endpoint $u=1$ is incompatible with
\eqref{eq:second}. Indeed,
\eqref{eq:C} and \eqref{eq:error-definition}
give
\[
C_{n-1}(1)-\Err_{n-1}(1)=1,
\]
whereas, for every $h>0$,
\[
\begin{aligned}
n\int_0^1
[\mathcal B_1(t)-ht(1-t)]_+^{(n-1)/2}
\dd t
&=
n\int_0^1
[(1-t)^2-ht(1-t)]_+^{(n-1)/2}
\dd t\\
&<
n\int_0^1(1-t)^{n-1}\dd t=1.
\end{aligned}
\]
The strict inequality holds because
$ht(1-t)>0$ for $0<t<1$.
Thus \eqref{eq:second} cannot hold
when $u=1$.

Moreover, $\tau_2=1$, so
\eqref{eq:polar-bounds} forces $u=1$
when $n=3$. Consequently,
\eqref{eq:polar} and \eqref{eq:second}
cannot hold simultaneously for $n=3$.
For the remainder of the argument,
we may therefore assume
$n\ge4$ and $0<u<1$.

We next derive a simpler upper bound
for the integral in \eqref{eq:second}.
Choose $\kappa>0$ such that
\[
h_0:=\frac{\kappa\eps}{u}<h,
\qquad
\frac{\kappa}{u}\le4,
\]
and define
\begin{equation}\label{eq:Q0-definition}
Q_0(t):=\mathcal B_u(t)-h_0t(1-t).
\end{equation}
For convenience, write
\begin{equation}\label{eq:Q0-parameters}
D:=u+h_0,\qquad
b:=u+\frac{h_0}{2},\qquad
m_0:=\eps-\frac{h_0^2}{4D}.
\end{equation}
Using \eqref{eq:Q0-definition} and
\eqref{eq:Q0-parameters}, we obtain
\begin{equation}\label{eq:Q0-square}
Q_0(t)
=1-2bt+Dt^2
=m_0+D\left(t-\frac{b}{D}\right)^2.
\end{equation}

Since $h_0\le4\eps$ and
$D=u+h_0\ge h_0>0$, we have
\[
m_0
=\eps-\frac{h_0^2}{4D}
\ge\eps-\frac{h_0}{4}\ge0.
\]
It follows from \eqref{eq:Q0-square}
that $Q_0(t)\ge0$ for every real $t$.\footnote{The parabola $Q_0$
has minimum height $m_0$ at $t=b/D$, with $D=Q_0''/2$ and
$b=-Q_0'(0)/2$. Subtracting $h_0t(1-t)$ from $\mathcal B_u(t)$ lowers
its interior while fixing both endpoint values.}

For integers $m\ge0$, define
\begin{equation}\label{eq:quadratic}
R_m:=(m+1)\int_0^1Q_0(t)^{m/2}\dd t.
\end{equation}
Since $h_0<h$ and $t(1-t)\ge0$ on
$[0,1]$, equation~\eqref{eq:Q0-definition}
gives
\[
\mathcal B_u(t)-ht(1-t)
\le Q_0(t).
\]
Together with $Q_0(t)\ge0$, this implies
\[
[\mathcal B_u(t)-ht(1-t)]_+
\le Q_0(t),
\qquad 0\le t\le1.
\]
Combining \eqref{eq:second} and
\eqref{eq:quadratic}, we conclude that
\begin{equation}\label{eq:quadratic-comparison}
C_{n-1}(u)-\Err_{n-1}(u)
\le R_{n-1}.
\end{equation}

We now verify that the choice
$
\kappa={6}/{(n-1)}$
satisfies both conditions imposed above.
By \eqref{eq:polar-bounds},
\[
h_0=\frac{6\eps}{(n-1)u}<h.
\]
To verify $\kappa/u\le4$, we use the
elementary numerical bounds
\begin{equation}\label{eq:e-bounds}
\frac52<\ee<
1+1+\frac12+\frac16+
\frac{1/24}{1-1/5}
=\frac{87}{32}<\frac{25}{9},
\qquad
(87/32)^3<(5/2)^4.
\end{equation}
The bounds for $\ee$ follow from its
exponential series, while the final
comparison is a rational inequality.
In particular, \eqref{eq:e-bounds}
implies $\ee^{3/4}<5/2$, and hence
$2\log(5/2)>3/2$.

Since $n\ge4$, the first inequality
in \eqref{eq:polar-bounds}, together
with $\ee^r>1+r$ for $r>0$, gives
\[
\begin{aligned}
(n-1)u
&\ge(n-1)
\left(
\left(\frac{n+1}{2}\right)^{2/(n-1)}
-1
\right)\\
&>2\log\frac{n+1}{2}\\
&\ge2\log\frac52>\frac32.
\end{aligned}
\]
Consequently,
\[
\frac{\kappa}{u}
=\frac{6}{(n-1)u}<4.
\]
Thus $\kappa=6/(n-1)$ is always an
admissible choice under
\eqref{eq:polar} for $n\ge4$ and
$0<u<1$.

The finite verification in
Appendix~\ref{app:verification}
also permits larger values of $\kappa$.
For a proposed value of $\kappa$,
let $h_0=\kappa\eps/u$ and suppose
$0<h_0\le1$.
Direct integration gives
\begin{equation}\label{eq:exact-polar}
P_{n-1}(u,h_0)
=
\frac{
2\bigl(
(1+h_0u)^{(n+1)/2}
-(1-h_0)^{(n+1)/2}
\bigr)
}{
(n+1)h_0(1+u)
}.
\end{equation}
We claim that
\begin{equation}\label{eq:polar-certification}
P_{n-1}(u,h_0)<1
\quad\Longrightarrow\quad h_0<h,
\end{equation}
provided that \eqref{eq:polar} holds.

Indeed, suppose instead that $h\le h_0$.
The convexity of $P_{n-1}(u,\cdot)$
established in the proof of
\eqref{eq:polar-chord} gives
\[
\begin{aligned}
P_{n-1}(u,h)
&\le
\left(1-\frac{h}{h_0}\right)
P_{n-1}(u,0)
+\frac{h}{h_0}P_{n-1}(u,h_0)\\
&=
1-\frac{h}{h_0}
+\frac{h}{h_0}P_{n-1}(u,h_0)
<1,
\end{aligned}
\]
contrary to \eqref{eq:polar}.
This proves \eqref{eq:polar-certification}.

Therefore, for any proposed
$\kappa>0$ satisfying
$\kappa/u\le4$ and $0<h_0\le1$,
the condition
$P_{n-1}(u,h_0)<1$ in
\eqref{eq:polar-certification}
provides a sufficient test for
$h_0<h$. In particular, no monotonicity
assumption on $P_{n-1}(u,h)$ as a
function of $h$ is required.

\section{Bounds for the centred product}\label{sec:error}

We derive two bounds for the centred-product error
defined in \eqref{eq:error-definition}.
The coefficient bound developed below is particularly
useful for small degrees or when $u$ is close to one.
A second bound, based on an integral estimate, will
be established in the next subsection.

\subsection{Coefficient bounds}

We begin with classical Schoenberg's inequality,
which was proved independently by Pereira
\cite{Per03} and Malamud \cite{Mal05}.
We shall use the following centred form.

\begin{lem}[Schoenberg's inequality]
Let $p$ be a complex polynomial of degree $n\ge2$, with zeros
$z_1,\ldots,z_n$ and critical points
$\zeta_1,\ldots,\zeta_{n-1}$, all counted with multiplicity.
If $\sum_{j=1}^n z_j=0$, then
\begin{equation}\label{eq:schoenberg-variance}
\sum_{k=1}^{n-1}|\zeta_k|^2
\le\frac{n-2}{n}\sum_{j=1}^n|z_j|^2,
\end{equation}
where equality holds if and only if all zeros are collinear in the complex plane.
\label{lem:schoenberg}
\end{lem}

For $\boldsymbol\eta=(\eta_1,\ldots,\eta_{n-1})
\in\mathbb C^{n-1}$, let
\[
e_k(\boldsymbol\eta)
=\sum_{1\le j_1<\cdots<j_k\le n-1}
\eta_{j_1}\cdots\eta_{j_k},
\qquad 1\le k\le n-1,
\]
and set $e_0(\boldsymbol\eta)=1$.

We shall bound the elementary symmetric functions
using differentiation, Cauchy's coefficient
estimate, and Newton's identities.
For the bound obtained from Newton's identities,
define $b_{n-1,0}=1$, $b_{n-1,1}=0$, and
\begin{equation}\label{eq:coefficients}
b_{n-1,k}
:=\frac{n-1}{k}
\sum_{j=2}^k
(n-2)^{\frac{j-2}{2}}b_{n-1,k-j},
\qquad 2\le k\le n-1.
\end{equation}

For $2\le k<n-1$, define
\begin{equation}\label{eq:coefficient-majorant}
\begin{aligned}
a_{n-1,k}:=\min\Biggl\{
\binom{n-1}{k}
\left(\frac{k-1}{n-2}\right)^{\frac{k}{2}},
\left(
\frac{(n-1)^{n-1}}
{k^k(n-1-k)^{n-1-k}}
\right)^{\frac12},
\quad b_{n-1,k}
\Biggr\},
\end{aligned}
\end{equation}
and set
$
a_{n-1,n-1}:=\min\{1,b_{n-1,n-1}\}.
$

\begin{lem}[Bounds for elementary symmetric functions]
\label{lem:coefficients}
Let $n\ge3$, and suppose that
$\eta_1,\ldots,\eta_{n-1}\in\mathbb C$
satisfy
\[
\sum_{j=1}^{n-1}\eta_j=0,
\qquad
\frac1{n-1}\sum_{j=1}^{n-1}|\eta_j|^2
\le V,
\]
where $V\ge0$.\footnote{Since the points $\eta_j$ have centroid
zero, $V$ bounds their mean squared distance from their centroid;
$\sqrt V$ bounds their root mean square radius, not necessarily the
distance of each individual point from the origin.}
Then, for $2\le k\le n-1$,
\begin{equation}\label{eq:elementary-symmetric-bound}
|e_k(\boldsymbol\eta)|
\le a_{n-1,k}V^{\frac{k}{2}}.
\end{equation}
\end{lem}

\begin{proof}
If $V=0$, then every $\eta_j$ vanishes,
and \eqref{eq:elementary-symmetric-bound}
is immediate. We may therefore assume $V>0$.

We first apply Lemma~\ref{lem:schoenberg}.
Consider the monic polynomial
$
F(z)=\prod_{j=1}^{n-1}(z-\eta_j).
$
Fix an integer $k$ with $2\le k\le n-1$.
The normalized derivative
$
\frac{k!}{(n-1)!}F^{(n-1-k)}(z)
$
is a monic polynomial of degree $k$.
Let $\xi_1,\ldots,\xi_k$ denote its zeros,
counted with multiplicity.

Differentiation preserves the centroid
of the zeros. Since
$\sum_{j=1}^{n-1}\eta_j=0$,
the zeros of every derivative of $F$
also have centroid zero.
Applying \eqref{eq:schoenberg-variance}
successively to the derivatives of $F$
gives
\begin{equation}\label{eq:iterated-variance}
\begin{aligned}
\frac1k\sum_{j=1}^k|\xi_j|^2
&\le
\frac{k-1}{n-2}
\frac1{n-1}\sum_{j=1}^{n-1}|\eta_j|^2\\
&\le\frac{k-1}{n-2}V.
\end{aligned}
\end{equation}

The constant coefficient of the normalized
derivative is
$
(-1)^k\frac{e_k(\boldsymbol\eta)}
{\binom{n-1}{k}}.
$
By Vieta's formulas,
\[
\frac{|e_k(\boldsymbol\eta)|}
{\binom{n-1}{k}}
=\prod_{j=1}^k|\xi_j|.
\]
The arithmetic--geometric mean inequality,
together with \eqref{eq:iterated-variance},
therefore yields
\begin{equation}\label{eq:coeff-schoenberg}
\begin{aligned}
|e_k(\boldsymbol\eta)|
&\le\binom{n-1}{k}
\left(\frac1k\sum_{j=1}^k|\xi_j|^2
\right)^{\frac{k}{2}}\\
&\le\binom{n-1}{k}
\left(\frac{k-1}{n-2}\right)^{\frac{k}{2}}
V^{\frac{k}{2}}.
\end{aligned}
\end{equation}

We next derive a bound from Cauchy's
coefficient estimate. Since
$\sum_{j=1}^{n-1}\eta_j=0$,
the arithmetic--geometric mean inequality
gives, for every $z\in\mathbb C$,
\begin{equation}\label{eq:polynomial-product-bound}
\begin{aligned}
|F(z)|
&\le
\left(
\frac1{n-1}\sum_{j=1}^{n-1}
|z-\eta_j|^2
\right)^{\frac{n-1}{2}}\\
&=
\left(
|z|^2+
\frac1{n-1}\sum_{j=1}^{n-1}|\eta_j|^2
\right)^{\frac{n-1}{2}}\\
&\le(|z|^2+V)^{\frac{n-1}{2}}.
\end{aligned}
\end{equation}

For $2\le k<n-1$, the coefficient
of $z^{n-1-k}$ in $F$ is
$(-1)^ke_k(\boldsymbol\eta)$.
Applying Cauchy's coefficient estimate
on the circle $|z|=r$ and using
\eqref{eq:polynomial-product-bound}, we obtain
\[
|e_k(\boldsymbol\eta)|
\le\frac{(r^2+V)^{\frac{n-1}{2}}}
{r^{n-1-k}},
\qquad r>0.
\]
The right-hand side is minimized when
$
r^2=(n-1-k)V/{k}.
$
Substituting this value of $r$ gives
\begin{equation}\label{eq:coeff-cauchy}
|e_k(\boldsymbol\eta)|
\le
\left(
\frac{(n-1)^{n-1}}
{k^k(n-1-k)^{n-1-k}}
\right)^{\frac12}
V^{\frac{k}{2}},
\qquad 2\le k<n-1.
\end{equation}
For $k=n-1$, applying the
arithmetic--geometric mean inequality
directly to $|\eta_1|^2,\ldots,
|\eta_{n-1}|^2$ gives
\[
|e_{n-1}(\boldsymbol\eta)|
=\prod_{j=1}^{n-1}|\eta_j|
\le V^{\frac{n-1}{2}}.
\]

Finally, we use Newton's identities.
The condition $\sum_{j=1}^{n-1}\eta_j=0$
implies
\[
\eta_i=-\sum_{\substack{1\le j\le n-1\\j\ne i}}
\eta_j.
\]
By the Cauchy--Schwarz inequality,
\[
|\eta_i|^2
\le(n-2)\sum_{j\ne i}|\eta_j|^2
\le(n-2)\bigl((n-1)V-|\eta_i|^2\bigr).
\]
Rearranging gives
\[
|\eta_i|^2\le(n-2)V,
\qquad 1\le i\le n-1.
\]
Consequently, for every integer
$2\le j\le n-1$,
\begin{equation}\label{eq:power-sum-bound}
\begin{aligned}
\left|\sum_{i=1}^{n-1}\eta_i^j\right|
&\le\sum_{i=1}^{n-1}|\eta_i|^j\\
&\le
\left(\max_{1\le i\le n-1}|\eta_i|
\right)^{j-2}
\sum_{i=1}^{n-1}|\eta_i|^2\\
&\le(n-1)(n-2)^{\frac{j-2}{2}}
V^{\frac j2}.
\end{aligned}
\end{equation}

Since $e_1(\boldsymbol\eta)=0$,
Newton's identities take the form
\begin{equation}\label{eq:newton-centred}
k e_k(\boldsymbol\eta)
=
\sum_{j=2}^k(-1)^{j-1}
e_{k-j}(\boldsymbol\eta)
\sum_{i=1}^{n-1}\eta_i^j,
\qquad 2\le k\le n-1.
\end{equation}
Starting from
$e_0(\boldsymbol\eta)=1$ and
$e_1(\boldsymbol\eta)=0$,
induction on $k$, using
\eqref{eq:power-sum-bound},
\eqref{eq:newton-centred}, and
\eqref{eq:coefficients}, gives
\begin{equation}\label{eq:coeff-newton}
|e_k(\boldsymbol\eta)|
\le b_{n-1,k}V^{\frac{k}{2}},
\qquad 2\le k\le n-1.
\end{equation}

Combining \eqref{eq:coeff-schoenberg},
\eqref{eq:coeff-cauchy}, and
\eqref{eq:coeff-newton} with the definition
of $a_{n-1,k}$ in
\eqref{eq:coefficient-majorant}, and using
the direct product bound when $k=n-1$,
proves \eqref{eq:elementary-symmetric-bound}.
\end{proof}

We now apply Lemma~\ref{lem:coefficients}
to the error term in \eqref{eq:error-definition}.
For any admissible vector
$\boldsymbol\eta$ in that definition,
the condition
$e_1(\boldsymbol\eta)=0$ gives
\begin{equation}\label{eq:centered-product-expansion}
\prod_{j=1}^{n-1}(1-ut+t\eta_j)
-(1-ut)^{n-1}=
\sum_{k=2}^{n-1}
e_k(\boldsymbol\eta)t^k
(1-ut)^{n-1-k}.
\end{equation}
Since $0<u\le1$ and $0\le t\le1$,
we have $1-ut\ge0$.
By \eqref{eq:elementary-symmetric-bound}
with $V=u\eps$,
\eqref{eq:centered-product-expansion} yields
\begin{equation}\label{eq:error-coefficient-integral}
\Err_{n-1}(u)
\le n\sum_{k=2}^{n-1}
a_{n-1,k}(u\eps)^{\frac{k}{2}}
\int_0^1 t^k(1-ut)^{n-1-k}\dd t.
\end{equation}

The integrals in
\eqref{eq:error-coefficient-integral}
can be evaluated explicitly.
Using
$
1-ut=(1-t)+\eps t
$
and the binomial theorem, we obtain
\[
\int_0^1t^k(1-ut)^{n-1-k}\dd t
=
\sum_{j=0}^{n-1-k}
\binom{n-1-k}{j}\eps^j
\int_0^1
t^{k+j}(1-t)^{n-1-k-j}\dd t.
\]
For nonnegative integers $k,j$ with
$k+j\le n-1$, the beta integral gives
\[
\int_0^1
t^{k+j}(1-t)^{n-1-k-j}\dd t
=
\frac{(k+j)!(n-1-k-j)!}{n!}.
\]
Substituting this evaluation into the
binomial expansion yields
\begin{equation}\label{eq:beta-integral}
n\int_0^1t^k(1-ut)^{n-1-k}\dd t
=
\frac1{\binom{n-1}{k}}
\sum_{j=0}^{n-1-k}
\binom{k+j}{k}\eps^j.
\end{equation}

Combining \eqref{eq:error-coefficient-integral}
and \eqref{eq:beta-integral}, we obtain
\begin{equation}\label{eq:error-finite}
\Err_{n-1}(u)
\le
\sum_{k=2}^{n-1}
\frac{a_{n-1,k}}{\binom{n-1}{k}}
(u\eps)^{\frac{k}{2}}
\sum_{j=0}^{n-1-k}
\binom{k+j}{k}\eps^j.
\end{equation}

A simpler bound follows when
$0<u<1$ and $\eps/u<1$.
Since $0<\eps<1$, the binomial series gives
\begin{equation}\label{eq:binomial-series-bound}
\begin{aligned}
\sum_{j=0}^{n-1-k}
\binom{k+j}{k}\eps^j
&\le
\sum_{j=0}^{\infty}
\binom{k+j}{k}\eps^j\\
&=(1-\eps)^{-k-1}
=u^{-k-1}.
\end{aligned}
\end{equation}
Applying \eqref{eq:coeff-schoenberg}
and \eqref{eq:binomial-series-bound}
to \eqref{eq:error-finite}, and dividing
by $\eps>0$, we obtain
\begin{equation}\label{eq:error-series}
\frac{\Err_{n-1}(u)}{\eps}
\le\frac1{u^2}
\sum_{k=2}^{n-1}
\left(\frac{k-1}{n-2}\right)^{\frac{k}{2}}
\left(\frac{\eps}{u}\right)^{\frac{k-2}{2}}.
\end{equation}

Fix an integer $k_*$ with
$2\le k_*\le n-1$.
The term corresponding to $k=2$
in \eqref{eq:error-series} equals
$1/(n-2)$.
For $k_*<k\le n-1$, we have
\[
0<
\left(\frac{k-1}{n-2}\right)^{\frac{k}{2}}
\le1.
\]
Consequently, when $k_*<n-1$,
the tail of the sum in
\eqref{eq:error-series} satisfies
\begin{equation}\label{eq:error-geometric-tail}
\sum_{k=k_*+1}^{n-1}
\left(\frac{k-1}{n-2}\right)^{\frac{k}{2}}
\left(\frac{\eps}{u}\right)^{\frac{k-2}{2}}\le
\sum_{k=k_*+1}^{\infty}
\left(\frac{\eps}{u}\right)^{\frac{k-2}{2}}=
\frac{
\left(\frac{\eps}{u}\right)^{\frac{k_*-1}{2}}
}{
1-\sqrt{\frac{\eps}{u}}
}.
\end{equation}

Combining \eqref{eq:error-series}
and \eqref{eq:error-geometric-tail}
gives
\begin{equation}\label{eq:error-geometric}
\frac{\Err_{n-1}(u)}{\eps}
\le\frac1{u^2}\Biggl(\frac1{n-2}
+\sum_{k=3}^{k_*}
\left(\frac{k-1}{n-2}\right)^{\frac{k}{2}}
\left(\frac{\eps}{u}\right)^{\frac{k-2}{2}}+
\frac{
\left(\frac{\eps}{u}\right)^{\frac{k_*-1}{2}}
}{
1-\sqrt{\frac{\eps}{u}}
}
\Biggr).
\end{equation}
When $k_*=n-1$, the sum has no tail,
and the final term in
\eqref{eq:error-geometric} is omitted.

\subsection{An integral bound}

We now derive an integral bound for the error
defined in \eqref{eq:error-definition}.
The main ingredient is a general form of
Zhang's centred product estimate
\cite[Lemma~5.1]{Zha26b}.

\begin{lem}[A centred product estimate]
\label{lem:interpolation}
Let $n\ge3$, $d\in\mathbb C\setminus\{0\}$,
and $\eta_1,\ldots,\eta_{n-1}\in\mathbb C$
satisfy
\[
\sum_{j=1}^{n-1}\eta_j=0,
\qquad
\frac1{n-1}\sum_{j=1}^{n-1}|\eta_j|^2\le V,
\]
where $V\ge0$. Then
\begin{equation}\label{eq:interpolation}
\left|
\prod_{j=1}^{n-1}(d+\eta_j)-d^{n-1}
\right|
\le
\frac{5(n-1)}{6}\frac{V}{|d|}
(|d|^2+V)^{\frac{n-2}{2}}.
\end{equation}
\end{lem}

\begin{proof}
For $0\le r\le1$, define
$
G(r)=\prod_{j=1}^{n-1}(d+r\eta_j).
$\footnote{The points $d+r\eta_j$ move along the segments from $d$ to
$d+\eta_j$. Their centroid remains $d$, and their mean squared distance
from it is at most $r^2V$. Thus varying $r$ contracts or expands the
configuration about its fixed centroid.}
Since $\sum_{j=1}^{n-1}\eta_j=0$,
we have
\[
0=\left(\sum_{j=1}^{n-1}\eta_j\right)G(r)=\sum_{j=1}^{n-1}
\eta_j(d+r\eta_j)
\prod_{k\ne j}(d+r\eta_k).
\]
Differentiating $G$ and rearranging
this identity gives
\begin{equation}\label{eq:centred-derivative}
dG'(r)
=-r\sum_{j=1}^{n-1}
\eta_j^2\prod_{k\ne j}(d+r\eta_k).
\end{equation}

We next estimate the products in
\eqref{eq:centred-derivative}.
Since the $\eta_j$ have mean zero,
\begin{equation}\label{eq:centred-mean-square}
\begin{aligned}
\sum_{k=1}^{n-1}|d+r\eta_k|^2
&=(n-1)|d|^2
+r^2\sum_{k=1}^{n-1}|\eta_k|^2\\
&\le(n-1)(|d|^2+V),
\end{aligned}
\end{equation}
where the inequality uses $0\le r\le1$.

For each $1\le j\le n-1$, the
arithmetic--geometric mean inequality
and \eqref{eq:centred-mean-square} yield
\begin{equation}\label{eq:centred-omitted-product}
\begin{aligned}
\left|\prod_{k\ne j}(d+r\eta_k)\right|
&\le
\left(
\frac1{n-2}\sum_{k\ne j}
|d+r\eta_k|^2
\right)^{\frac{n-2}{2}}\\
&\le
\left(\frac{n-1}{n-2}\right)^{\frac{n-2}{2}}
(|d|^2+V)^{\frac{n-2}{2}}.
\end{aligned}
\end{equation}

Combining \eqref{eq:centred-derivative}
and \eqref{eq:centred-omitted-product},
and using
$\sum_{j=1}^{n-1}|\eta_j|^2\le(n-1)V$,
we obtain
\[
|G'(r)|
\le
\frac{(n-1)rV}{|d|}
\left(\frac{n-1}{n-2}\right)^{\frac{n-2}{2}}
(|d|^2+V)^{\frac{n-2}{2}}.
\]
Since $G(0)=d^{n-1}$ and
$G(1)=\prod_{j=1}^{n-1}(d+\eta_j)$,
integration over $[0,1]$ gives
\begin{equation}\label{eq:centred-integrated-bound}
\begin{aligned}
\left|
\prod_{j=1}^{n-1}(d+\eta_j)-d^{n-1}
\right|
&\le
\int_0^1|G'(r)|\dd r\\
&\le
\frac{n-1}{2}\frac{V}{|d|}
\left(\frac{n-1}{n-2}\right)^{\frac{n-2}{2}}
(|d|^2+V)^{\frac{n-2}{2}}.
\end{aligned}
\end{equation}
Finally,
\[
\left(\frac{n-1}{n-2}\right)^{\frac{n-2}{2}}
<\sqrt{\ee}<\frac53,
\]
where the last inequality follows from
\eqref{eq:e-bounds}.
Substitution into \eqref{eq:centred-integrated-bound}
proves \eqref{eq:interpolation}.
\end{proof}

We apply Lemma~\ref{lem:interpolation}
to the centred-product error
$\Err_{n-1}(u)$. The case $u=1$ is
immediate from \eqref{eq:error-definition},
since $\Err_{n-1}(1)=0$.
We therefore assume $0<u<1$.

Fix an arbitrary vector
$\boldsymbol\eta=(\eta_1,\ldots,\eta_{n-1})$
satisfying the constraints in
\eqref{eq:error-definition}.
For $0\le t\le1$, observe that
\begin{equation}\label{eq:B-identity}
\mathcal B_u(t)
=(1-ut)^2+u\eps t^2.
\end{equation}
Since $1-ut\ge1-u=\eps>0$,
we may apply \eqref{eq:interpolation}
with $d=1-ut$, with $t\eta_j$
in place of $\eta_j$, and with
$V=u\eps t^2$.
Using \eqref{eq:B-identity}, we obtain
\begin{equation}\label{eq:error-pointwise-interpolation}
\left|
\prod_{j=1}^{n-1}(1-ut+t\eta_j)
-(1-ut)^{n-1}
\right|\le
\frac{5(n-1)}6
u\eps\frac{t^2}{1-ut}
\mathcal B_u(t)^{\frac{n-2}{2}}.
\end{equation}

A second pointwise bound follows from
the arithmetic--geometric mean inequality.
Since $\sum_{j=1}^{n-1}\eta_j=0$,
we have
\[
\begin{aligned}
\frac1{n-1}\sum_{j=1}^{n-1}
|1-ut+t\eta_j|^2
&=(1-ut)^2+
\frac{t^2}{n-1}\sum_{j=1}^{n-1}|\eta_j|^2\\
&\le(1-ut)^2+u\eps t^2
=\mathcal B_u(t),
\end{aligned}
\]
where the last equality uses
\eqref{eq:B-identity}.
Consequently,
\begin{equation}\label{eq:error-pointwise-product}
\left|
\prod_{j=1}^{n-1}(1-ut+t\eta_j)
\right|
\le\mathcal B_u(t)^{\frac{n-1}{2}}.
\end{equation}
The triangle inequality and
\eqref{eq:error-pointwise-product} give
\begin{equation}\label{eq:error-pointwise-triangle}
\left|
\prod_{j=1}^{n-1}(1-ut+t\eta_j)
-(1-ut)^{n-1}
\right|\le
\mathcal B_u(t)^{\frac{n-1}{2}}
+(1-ut)^{n-1}.
\end{equation}

For $0\le t\le1$, define
\begin{equation}\label{eq:OT}
\begin{aligned}
O_{n-1}(t)
&:=\frac56n(n-1)u\eps
\frac{t^2}{1-ut}
\mathcal B_u(t)^{\frac{n-2}{2}},\\
T_{n-1}(t)
&:=n\left(
\mathcal B_u(t)^{\frac{n-1}{2}}
+(1-ut)^{n-1}
\right).
\end{aligned}
\end{equation}
By \eqref{eq:error-pointwise-interpolation},
\eqref{eq:error-pointwise-triangle}, and
\eqref{eq:OT}, every admissible vector
$\boldsymbol\eta$ satisfies
\[
n\left|
\prod_{j=1}^{n-1}(1-ut+t\eta_j)
-(1-ut)^{n-1}
\right|\le
\min\{O_{n-1}(t),T_{n-1}(t)\},
\qquad 0\le t\le1.
\]
The right-hand side is independent of
$\boldsymbol\eta$. Integrating over
$[0,1]$ and taking the maximum in
\eqref{eq:error-definition}, we obtain
\begin{equation}\label{eq:error-integral}
\Err_{n-1}(u)
\le\int_0^1
\min\{O_{n-1}(t),T_{n-1}(t)\}\dd t.
\end{equation}

We next derive upper bounds for the
integrals in \eqref{eq:error-integral}.
Throughout the remaining argument, assume
$n\ge4$ and $0<u<1$.

Since
\[
\mathcal B_u(t)
=1-ut(2-t)\le1-ut,
\qquad 0\le t\le1,
\]
and $\mathcal B_u(t)\ge\eps>0$,
the definition of $O_{n-1}(t)$ in
\eqref{eq:OT} gives
\begin{equation}\label{eq:O-majorant}
\frac{O_{n-1}(t)}{\eps}
\le\frac56n(n-1)u\,t^2
\mathcal B_u(t)^{\frac{n-4}{2}}.
\end{equation}
Here the exponent $(n-4)/2$ is
nonnegative because $n\ge4$.

We shall use the elementary beta integral
\begin{equation}\label{eq:quadratic-beta-integral}
\int_0^{1/c}
t^2(1-ct)^{\frac{n-4}{2}}\dd t
=
\frac{16}{c^3(n-2)n(n+2)},
\qquad c>0.
\end{equation}
Indeed, substituting $v=ct$ gives
\[
\int_0^{1/c}
t^2(1-ct)^{\frac{n-4}{2}}\dd t=
\frac1{c^3}
\int_0^1v^2(1-v)^{\frac{n-4}{2}}\dd v=
\frac1{c^3}
\frac{2}{
\frac{n-2}{2}\cdot
\frac n2\cdot
\frac{n+2}{2}},
\]
which proves \eqref{eq:quadratic-beta-integral}.

For the integral over $[0,1]$, use
$\mathcal B_u(t)\le1-ut$ in
\eqref{eq:O-majorant}.
Since $1/u>1$, the integral can be
extended to $[0,1/u]$.
Combining \eqref{eq:error-integral},
\eqref{eq:O-majorant}, and
\eqref{eq:quadratic-beta-integral}
with $c=u$, we obtain
\begin{align}
\frac{\Err_{n-1}(u)}{\eps}
&\le\frac1\eps\int_0^1O_{n-1}(t)\dd t
\notag\\
&\le\frac56n(n-1)u
\int_0^{1/u}
t^2(1-ut)^{\frac{n-4}{2}}\dd t
\notag\\
&=\frac{40(n-1)}
{3(n-2)(n+2)u^2}.
\label{eq:error-whole}
\end{align}

Now fix a splitting point $d\in(0,1)$.
For $0\le t\le d$, we have
\[
\mathcal B_u(t)
=1-ut(2-t)
\le1-(2-d)ut.
\]
Moreover, $ud(2-d)<1$, so
\[
d<\frac1{(2-d)u}.
\]
Applying \eqref{eq:O-majorant}
and extending the integral to
$[0,1/((2-d)u)]$, we obtain
\begin{align}
\frac1\eps\int_0^dO_{n-1}(t)\dd t
&\le\frac56n(n-1)u
\int_0^{1/((2-d)u)}
t^2[1-(2-d)ut]^{\frac{n-4}{2}}
\dd t
\notag\\
&=\frac{40(n-1)}
{3(n-2)(n+2)(2-d)^3u^2},
\label{eq:error-front}
\end{align}
where the equality follows from
\eqref{eq:quadratic-beta-integral}
with $c=(2-d)u$.

It remains to estimate the integral
over $[d,1]$.
Both $\mathcal B_u(t)$ and $1-ut$
are nonincreasing on $[0,1]$, since
\[
\mathcal B_u'(t)=-2u(1-t)\le0,
\qquad
\frac{\dd}{\dd t}(1-ut)=-u<0.
\]
Using $t^2\le1$ in \eqref{eq:O-majorant},
we obtain, for $d\le t\le1$,
\[
\frac{O_{n-1}(t)}{\eps}
\le
\frac56n(n-1)u
\mathcal B_u(d)^{\frac{n-4}{2}}.
\]
Similarly, the definition of $T_{n-1}(t)$
in \eqref{eq:OT} gives
\[
\frac{T_{n-1}(t)}{\eps}
\le
\frac n\eps
\left(
\mathcal B_u(d)^{\frac{n-1}{2}}
+(1-ud)^{n-1}
\right).
\]

Integrating the two pointwise bounds
over $[d,1]$, we conclude that
\begin{equation}\label{eq:error-tail}
\frac1\eps\int_d^1
\min\{O_{n-1}(t),T_{n-1}(t)\}\dd t\le\min\Biggl\{
\frac56n(n-1)u(1-d)
\mathcal B_u(d)^{\frac{n-4}{2}},
\frac{n(1-d)}{\eps}
\left(
\mathcal B_u(d)^{\frac{n-1}{2}}
+(1-ud)^{n-1}
\right)
\Biggr\}.
\end{equation}

\section{Quadratic integrals}\label{sec:moments}

\begin{figure}[tbp]
\centering
\includegraphics[width=\textwidth]{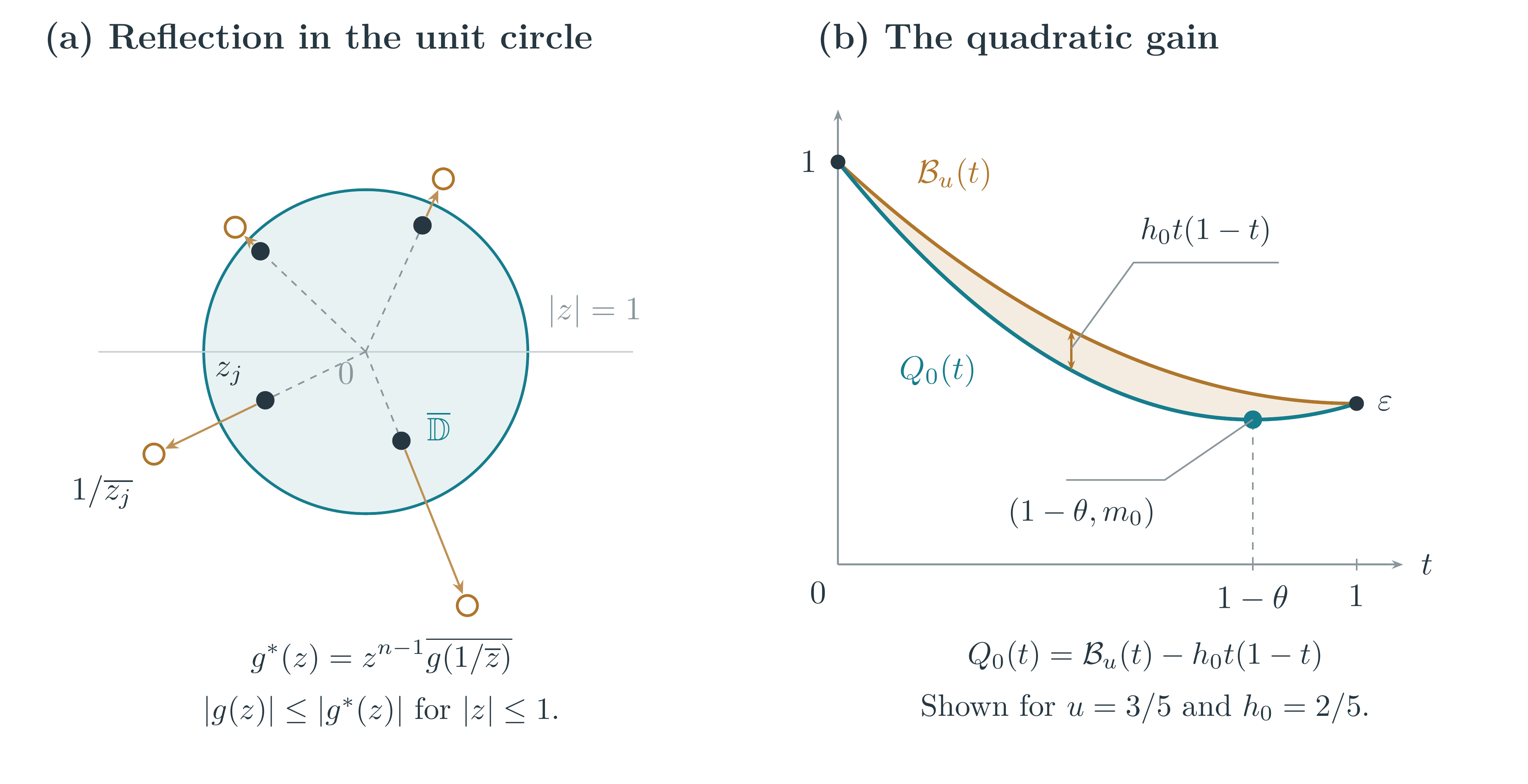}
\caption{Reflection and the quadratic comparison.
(a) Each nonzero zero $z_j$ of $g(z)=p(z)/(z-a)$ is
reflected to the zero $1/\overline{z_j}$ of $g^*$.
Filled and open circles represent the zeros of $g$ and
$g^*$, respectively.
(b) Subtracting $h_0t(1-t)$ from $\mathcal B_u(t)$
lowers the graph on $(0,1)$ and preserves both endpoint
values. The minimum of $Q_0$ occurs at $t=1-\theta$.
Here $u=3/5$, $h_0=2/5$, $\eps=2/5$,
$\theta=1/5$, and $m_0=9/25$.}
\label{fig:reflection-quadratic}
\end{figure}

Throughout this section, assume $0<u<1$ and choose
$\kappa>0$ such that
\[
h_0=\frac{\kappa\eps}{u}<h,
\qquad \frac{\kappa}{u}\le4,
\]
where $\eps=1-u$.
Let $Q_0,D,b,m_0$, and $R_m$ be defined by
\eqref{eq:Q0-definition}, \eqref{eq:Q0-parameters},
and \eqref{eq:quadratic}. In particular,
\eqref{eq:Q0-square} gives
\[
Q_0(t)=m_0+D\left(t-\frac{b}{D}\right)^2,
\qquad m_0\ge0.
\]
Set
\begin{equation}\label{eq:theta-definition}
\theta:=\frac{h_0}{2D},
\qquad
\frac{b}{D}=1-\theta,
\end{equation}
where the second identity follows from
\eqref{eq:Q0-parameters}.\footnote{The number $\theta$ is the horizontal
distance from the vertex of the parabola $Q_0$ to the right endpoint
$t=1$. It is the length of the right-hand interval when the integral is
split at that vertex.}
Since $0<h_0<D$, we have $0<\theta<1/2$. Figure~\ref{fig:reflection-quadratic} illustrates reflection
in the unit circle and the quadratic comparison
in~\eqref{eq:Q0-definition}.

\begin{lem}[A recurrence for quadratic integrals]
\label{lem:recurrence}
For every integer $m\ge2$,
\begin{equation}\label{eq:recurrence}
R_m
=1-\theta+\theta\eps^{\frac m2}
+\frac{m}{m-1}m_0R_{m-2}.
\end{equation}
\end{lem}

\begin{proof}
By \eqref{eq:Q0-square},
\[
Q_0'(t)=2D\left(t-\frac bD\right).
\]
Differentiating
$\left(t-b/D\right)Q_0(t)^{\frac m2}$
and using
\[
D\left(t-\frac bD\right)^2=Q_0(t)-m_0,
\]
we obtain
\begin{equation}\label{eq:quadratic-derivative}
\begin{aligned}
\frac{\dd}{\dd t}
\left[
\left(t-\frac bD\right)
Q_0(t)^{\frac m2}
\right]&=
Q_0(t)^{\frac m2}
+mD\left(t-\frac bD\right)^2
Q_0(t)^{\frac{m-2}{2}}\\
&=
(m+1)Q_0(t)^{\frac m2}
-mm_0Q_0(t)^{\frac{m-2}{2}}.
\end{aligned}
\end{equation}

By \eqref{eq:Q0-definition},
\[
Q_0(0)=1,
\qquad Q_0(1)=\eps.
\]
Consequently, \eqref{eq:theta-definition} gives
\[
\begin{aligned}
\left[
\left(t-\frac bD\right)
Q_0(t)^{\frac m2}
\right]_{t=0}^{t=1}
&=
\left(1-\frac bD\right)\eps^{\frac m2}
+\frac bD\\
&=
1-\theta+\theta\eps^{\frac m2}.
\end{aligned}
\]
Integrating \eqref{eq:quadratic-derivative}
over $[0,1]$ and using
\eqref{eq:quadratic}, we obtain
\[
R_m-\frac{m}{m-1}m_0R_{m-2}
=1-\theta+\theta\eps^{\frac m2}.
\]
Rearranging proves \eqref{eq:recurrence}.
\end{proof}

To bound $R_m$, introduce the auxiliary integral
\begin{equation}\label{eq:F}
\mathcal F_m(\alpha)
:=(m+1)\int_0^1
(1-\alpha+\alpha t^2)^{\frac m2}\dd t,
\qquad 0\le\alpha\le1,\quad m\in\mathbb Z_{\ge0}.
\end{equation}
Since $h_0>0$, equation
\eqref{eq:Q0-definition} gives
\[
Q_0(t)\le\mathcal B_u(t)
=\eps+u(1-t)^2,
\qquad 0\le t\le1.
\]
Applying the change of variable $t\mapsto1-t$
to \eqref{eq:quadratic}, we obtain
\begin{equation}\label{eq:F-comparison}
\begin{aligned}
R_m
&\le(m+1)\int_0^1
[\eps+u(1-t)^2]^{\frac m2}\dd t\\
&=\mathcal F_m(u).
\end{aligned}
\end{equation}

We establish three elementary bounds for
$\mathcal F_m(u)$. For every integer $m\ge2$,
\begin{equation}\label{eq:F-elementary}
\begin{aligned}
\mathcal F_m(u)&\le1+m\eps,\\
\mathcal F_m(u)&<\frac2u,\\
\mathcal F_m(u)&\le1+3\eps
\qquad\text{if }\eps\le\frac14.
\end{aligned}
\end{equation}

For the first inequality in
\eqref{eq:F-elementary}, write
\[
\mathcal F_m(u)
=(m+1)\int_0^1
[t^2+\eps(1-t^2)]^{\frac m2}\dd t.
\]
Since $m/2\ge1$, the integrand is
convex as a function of $\eps\in[0,1]$.
At $\eps=0$ and $\eps=1$, the integral
has values $1/(m+1)$ and $1$,
respectively. Convexity therefore gives
\[
\mathcal F_m(u)
\le(1-\eps)+\eps(m+1)=1+m\eps.
\]

For the second inequality in
\eqref{eq:F-elementary}, use $t^2\le t$
for $0\le t\le1$. Direct integration gives
\[
\begin{aligned}
\mathcal F_m(u)
&\le(m+1)\int_0^1
(1-u+ut)^{\frac m2}\dd t\\
&=\frac{2(m+1)}{u(m+2)}
\left(1-\eps^{\frac{m+2}{2}}\right)\\
&<\frac2u.
\end{aligned}
\]

For the third inequality, integration
by parts in \eqref{eq:F} gives
\begin{equation}\label{eq:F-recurrence}
\mathcal F_m(\alpha)
=1+\frac{m}{m-1}(1-\alpha)
\mathcal F_{m-2}(\alpha),
\qquad m\ge2.
\end{equation}
Indeed, differentiating
$t(1-\alpha+\alpha t^2)^{\frac m2}$
and integrating over $[0,1]$ gives
\[
1
=(m+1)\int_0^1
(1-\alpha+\alpha t^2)^{\frac m2}\dd t
-m(1-\alpha)\int_0^1
(1-\alpha+\alpha t^2)^{\frac{m-2}{2}}
\dd t,
\]
which is equivalent to
\eqref{eq:F-recurrence}.

We have
\[
\mathcal F_2(u)=1+2\eps,
\qquad
\mathcal F_1(u)\le2.
\]
Applying \eqref{eq:F-recurrence} with
$m=3$ gives
$\mathcal F_3(u)\le1+3\eps$.
For $m\ge4$, assume that
$\mathcal F_{m-2}(u)\le1+3\eps$.
If $\eps\le1/4$, then
\eqref{eq:F-recurrence} yields
\[
\begin{aligned}
\mathcal F_m(u)
&\le1+\frac{m}{m-1}
\eps(1+3\eps)\\
&\le1+\frac43\eps
\left(1+\frac34\right)\\
&\le1+3\eps.
\end{aligned}
\]
Induction proves the third inequality
in \eqref{eq:F-elementary}.

We next establish a sharper bound for
$\mathcal F_m(\alpha)$ when $m\alpha$
is large.

\begin{lem}[An upper bound for quadratic integrals]
\label{lem:moment}
Let $m\ge2$ be an integer and
$0<\alpha\le1$. If $m\alpha\ge90$, then
\begin{equation}\label{eq:moment}
\mathcal F_m(\alpha)
\le\frac1\alpha+
\frac{11}{10}\frac{1-\alpha}{m\alpha^2}.
\end{equation}
\end{lem}

\begin{proof}
\emph{Step 1. Estimating the coefficients
in the recurrence.}
Put $r=1-\alpha$ and
$k_*=\lfloor m/4\rfloor$.
For $0\le j\le\lfloor m/2\rfloor$,
define
\[
c_0=1,
\qquad
c_j=\prod_{\ell=0}^{j-1}
\frac{m-2\ell}{m-2\ell-1}
\quad(j\ge1).
\]
Iterating \eqref{eq:F-recurrence}
until the index reaches zero or one,
and using
\[
\mathcal F_0(\alpha)=1,
\qquad
\mathcal F_1(\alpha)\le2,
\]
we obtain
\begin{equation}\label{eq:moment-expansion}
\mathcal F_m(\alpha)\le
\begin{cases}
\displaystyle
\sum_{j=0}^{m/2}c_jr^j,
& m\text{ even},\\[8pt]
\displaystyle
\sum_{j=0}^{(m-3)/2}c_jr^j
+2c_{(m-1)/2}r^{\frac{m-1}{2}},
& m\text{ odd}.
\end{cases}
\end{equation}
All coefficients in \eqref{eq:moment-expansion}
are nonnegative.

The sum of the coefficients in the
polynomial on the right-hand side of
\eqref{eq:moment-expansion} is $m+1$.
Indeed, setting $r=1$ corresponds to
$\alpha=0$, where
$\mathcal F_m(0)=m+1$ and the terminal
values $\mathcal F_0(0)=1$ and
$\mathcal F_1(0)=2$ are exact.

We estimate $c_j$ for $0\le j\le k_*$.
Using $\log(1+x)\le x$ for $x\ge0$,
we obtain
\[
\begin{aligned}
\log c_j
&=\sum_{\ell=0}^{j-1}
\log\left(1+\frac1{m-2\ell-1}\right)\\
&\le\sum_{\ell=0}^{j-1}
\frac1{m-2\ell-1}\\
&\le\frac jm+
\frac{j^2}{m(m-2j+1)}.
\end{aligned}
\]
For the last inequality, write
\[
\frac1{m-2\ell-1}
=\frac1m+
\frac{2\ell+1}{m(m-2\ell-1)}
\]
and use
$m-2\ell-1\ge m-2j+1$.
Since $j\le m/4$, we have
$m-2j+1\ge m/2$, and hence
\[
\log c_j
\le\frac jm+2\left(\frac jm\right)^2
\le\frac38.
\]

For $0\le x\le3/8$, the exponential
series gives
\[
\ee^x
\le1+x+\frac{x^2}{2(1-x/3)}.
\]
Applying this inequality with
$x=j/m+2(j/m)^2$, and using
$x\le3j/(2m)$, we obtain
\[
\begin{aligned}
c_j
&\le1+\frac jm
+2\frac{j^2}{m^2}
+\frac{9j^2}{7m^2}\\
&=1+\frac jm+\frac{23j^2}{7m^2}\\
&\le1+\frac jm+4\frac{j^2}{m^2}.
\end{aligned}
\]
Thus
\begin{equation}\label{eq:moment-coefficient-bound}
c_j\le
1+\frac jm+4\frac{j^2}{m^2},
\qquad 0\le j\le k_*.
\end{equation}

\medskip
\noindent
\emph{Step 2. Summing the coefficients.}
Since the coefficients in
\eqref{eq:moment-expansion} are nonnegative
and their sum is $m+1$, we have
\[
\mathcal F_m(\alpha)
\le\sum_{j=0}^{k_*}c_jr^j
+(m+1)r^{k_*+1}.
\]
Applying \eqref{eq:moment-coefficient-bound}
and extending the finite sum to an
infinite series gives
\[
\begin{aligned}
\mathcal F_m(\alpha)
&\le
\sum_{j=0}^{\infty}
\left(1+\frac jm+4\frac{j^2}{m^2}\right)r^j
+(m+1)r^{k_*+1}\\
&=\frac1\alpha+
\frac{r}{m\alpha^2}
+\frac{4r(1+r)}{m^2\alpha^3}
+(m+1)r^{k_*+1},
\end{aligned}
\]
where the equality uses
\[
\sum_{j=0}^{\infty}r^j=\frac1\alpha,
\qquad
\sum_{j=0}^{\infty}jr^j=\frac r{\alpha^2},
\qquad
\sum_{j=0}^{\infty}j^2r^j
=\frac{r(1+r)}{\alpha^3}.
\]
Since $1+r\le2$, we conclude that
\begin{equation}\label{eq:moment-series}
\mathcal F_m(\alpha)
\le\frac1\alpha+
\frac{r}{m\alpha^2}
\left(1+\frac8{m\alpha}\right)
+(m+1)r^{k_*+1}.
\end{equation}

If $r=0$, then $\alpha=1$ and
$\mathcal F_m(1)=1$, so
\eqref{eq:moment} holds with equality.
We therefore assume $0<r<1$.

Since $k_*=\lfloor m/4\rfloor$,
we have $k_*\ge m/4-1$.
Moreover, $\log r\le-\alpha$, so
\[
r^{k_*}
\le\exp(-\alpha k_*)
\le\ee^\alpha\exp\left(-\frac{m\alpha}{4}\right).
\]
It follows that
\begin{equation}\label{eq:moment-tail-bound}
\begin{aligned}
\frac{(m+1)r^{k_*+1}}
{r/(m\alpha^2)}
&=(m+1)m\alpha^2r^{k_*}\\
&\le
2\ee(m\alpha)^2
\exp\left(-\frac{m\alpha}{4}\right)\\
&<
6(m\alpha)^2
\exp\left(-\frac{m\alpha}{4}\right)\\
&<\frac1{1000}.
\end{aligned}
\end{equation}
For the second inequality, we used
$m\ge90$ and $0<\alpha\le1$.
The third follows from $\ee<3$
in \eqref{eq:e-bounds}.
For the last inequality, observe that
$x^2\ee^{-\frac x4}$ is decreasing for
$x\ge8$. Since $m\alpha\ge90$,
\[
6(m\alpha)^2
\exp\left(-\frac{m\alpha}{4}\right)
\le6\cdot90^2\ee^{-\frac{90}{4}}
<\frac1{1000},
\]
where the final numerical comparison
follows from $\ee>5/2$
in \eqref{eq:e-bounds}.

Combining \eqref{eq:moment-series}
and \eqref{eq:moment-tail-bound}, we obtain
\[
\begin{aligned}
\mathcal F_m(\alpha)
&\le\frac1\alpha+
\frac{r}{m\alpha^2}
\left(1+\frac8{m\alpha}
+\frac1{1000}\right)\\
&\le\frac1\alpha+
\frac{r}{m\alpha^2}
\left(1+\frac8{90}
+\frac1{1000}\right)\\
&\le\frac1\alpha+
\frac{11}{10}\frac{r}{m\alpha^2}.
\end{aligned}
\]
Recalling that $r=1-\alpha$ proves
\eqref{eq:moment}.
\end{proof}

We now apply Lemma~\ref{lem:moment}
to the quadratic $Q_0$.
By \eqref{eq:Q0-square}, define
\begin{equation}\label{eq:vertex-alpha}
\alpha:=\frac{b^2}{D}=1-m_0.
\end{equation}
Since $b,D>0$ and $m_0\ge0$,
we have $0<\alpha\le1$.\footnote{Here
$\alpha=Q_0(0)-\min_{t\in\mathbb R}Q_0(t)$ is
the vertical drop from the left endpoint to the vertex. Rescaling the
horizontal interval from the vertex to the left endpoint gives the
quadratic $1-\alpha+\alpha t^2$ in \eqref{eq:F}.}

The vertex of $Q_0$ is located at
$t=b/D=1-\theta$, which lies in $(0,1)$
by \eqref{eq:theta-definition}.
We split the integral in
\eqref{eq:quadratic} at $t=b/D$.

For the integral over $[0,b/D]$,
make the substitution
\[
t=\frac bD(1-v),
\qquad 0\le v\le1.
\]
Equations \eqref{eq:Q0-square}
and \eqref{eq:vertex-alpha} give
\[
Q_0\left(\frac bD(1-v)\right)
=m_0+\frac{b^2}{D}v^2
=1-\alpha+\alpha v^2.
\]
Therefore, by \eqref{eq:F},
\[
\begin{aligned}
(m+1)\int_0^{b/D}
Q_0(t)^{\frac m2}\dd t
&=
(m+1)\frac bD
\int_0^1(1-\alpha+\alpha v^2)^{\frac m2}
\dd v\\
&=\frac bD\mathcal F_m(\alpha).
\end{aligned}
\]

On $[b/D,1]$, the function $Q_0$
is nondecreasing, because
\[
Q_0'(t)=2D\left(t-\frac bD\right)\ge0.
\]
Since $Q_0(1)=\eps$ and
$1-b/D=\theta$, we obtain
\[
(m+1)\int_{b/D}^1
Q_0(t)^{\frac m2}\dd t
\le(m+1)\theta\eps^{\frac m2}.
\]
Combining the integrals over
$[0,b/D]$ and $[b/D,1]$, we conclude that
\begin{equation}\label{eq:R-split}
R_m\le
\frac bD\mathcal F_m(\alpha)
+(m+1)\theta\eps^{\frac m2}.
\end{equation}

Suppose now that
$(n-1)\alpha\ge90$.
Applying Lemma~\ref{lem:moment}
with $m=n-1$ to \eqref{eq:R-split},
we obtain
\[
\begin{aligned}
R_{n-1}
&\le
\frac bD\left(
\frac1\alpha+
\frac{11}{10}
\frac{1-\alpha}{(n-1)\alpha^2}
\right)
+n\theta\eps^{\frac{n-1}{2}}.
\end{aligned}
\]
By \eqref{eq:vertex-alpha},
\[
\frac b{D\alpha}=\frac1b,
\qquad
\frac{b(1-\alpha)}{D\alpha^2}
=\frac{m_0D}{b^3}.
\]
Substitution gives
\begin{equation}\label{eq:right-moment}
R_{n-1}
\le
\frac1b+
\frac{11}{10}\frac{m_0D}{(n-1)b^3}
+n\theta\eps^{\frac{n-1}{2}},
\qquad (n-1)\alpha\ge90.
\end{equation}

\section{Exclusion of the necessary inequalities}
\label{sec:exclusion}

We prove that the two necessary inequalities
\eqref{eq:polar} and \eqref{eq:second} cannot hold
simultaneously. The argument is divided into two
parts. For $n\ge6001$, we use the analytic bounds
established in Sections~\ref{sec:error} and
\ref{sec:moments}. The degrees $4\le n\le6000$
are handled by the rigorous interval verification
described in Appendix~\ref{app:verification}.

\subsection{The large-degree estimate}

\begin{lem}[Exclusion in large degrees]
\label{lem:large}
For every integer $n\ge6001$, there are no
$0<u<1$ and $0<h<1$ satisfying both
\eqref{eq:polar} and \eqref{eq:second}.
\end{lem}

\begin{proof}
Suppose, to the contrary, that
$0<u<1$ and $0<h<1$ satisfy both
\eqref{eq:polar} and \eqref{eq:second}.
We derive incompatible bounds for
$C_{n-1}(u)-R_{n-1}$ and
$\Err_{n-1}(u)$.

\medskip
\noindent
\emph{Step 1. A lower bound for
$C_{n-1}(u)-R_{n-1}$.}
Choose
\[
\kappa=\frac6{n-1},
\qquad
h_0=\frac{6\eps}{(n-1)u}.
\]
The argument following
\eqref{eq:quadratic-comparison} shows that
this choice satisfies $h_0<h$ and
$\kappa/u<4$ whenever \eqref{eq:polar}
holds. Thus the definitions in
\eqref{eq:Q0-parameters},
\eqref{eq:theta-definition}, and
\eqref{eq:vertex-alpha} apply.

We first establish a lower bound for
$(n-1)u$. From \eqref{eq:polar-bounds}
and $\ee^x>1+x$ for $x>0$, we have
\[
\begin{aligned}
(n-1)u
&\ge(n-1)
\left[
\left(\frac{n+1}{2}\right)^{\frac2{n-1}}
-1
\right]\\
&>2\log\frac{n+1}{2}
\ge2\log3001.
\end{aligned}
\]
By \eqref{eq:e-bounds} and the rational
inequality $(87/32)^8<3001$, we have
$\ee^8<3001$. Consequently,
\begin{equation}\label{eq:large-u}
t:=(n-1)u>16,
\qquad
h_0=\frac{6\eps}{t}<\frac{3\eps}{8}.
\end{equation}

Recall from \eqref{eq:vertex-alpha} that
$\alpha=b^2/D$, where
$
b=u+{h_0}/{2}$,
and  $D=u+h_0=2b-u.
$
Direct calculation gives
\[
\alpha-u=\frac{h_0^2}{4D}\ge0.
\]
Moreover, $D\le2b$, so
$\alpha\ge b/2$.

Suppose first that $0<u\le1/2$.
Since $\eps\ge1/2$, we obtain
\[
\begin{aligned}
\alpha
&\ge\frac b2
=\frac u2+\frac{3\eps}{2(n-1)u}\\
&\ge\frac u2+\frac3{4(n-1)u}\\
&\ge\sqrt{\frac3{2(n-1)}},
\end{aligned}
\]
where the last inequality is the
arithmetic--geometric mean inequality.
Thus
\[
(n-1)\alpha
\ge\sqrt{\frac{3(n-1)}2}
\ge\sqrt{9000}>90.
\]
If $1/2\le u<1$, then
$\alpha\ge u\ge1/2$, and hence
$(n-1)\alpha\ge3000>90$.
We have therefore proved
\begin{equation}\label{eq:large-alpha}
(n-1)\alpha>90
\qquad(n\ge6001).
\end{equation}
By \eqref{eq:large-alpha}, Lemma~\ref{lem:moment} applies,
and \eqref{eq:right-moment} gives
\[
R_{n-1}
\le\frac1b+
\frac{11}{10}\frac{m_0D}{(n-1)b^3}
+n\theta\eps^{\frac{n-1}{2}}.
\]

Introduce the normalization
\begin{equation}\label{eq:K0}
K_0:=\frac{\eps}{(n-1)u^2}>0.
\end{equation}
We estimate the two error terms
in \eqref{eq:right-moment} after
division by $K_0$.

Since $\theta<1/2$ by
\eqref{eq:theta-definition} and
$1-u\le\ee^{-u}$, we obtain
\[
\begin{aligned}
\frac{n\theta\eps^{\frac{n-1}{2}}}{K_0}
&\le
\frac{n(n-1)u^2}{2}
\eps^{\frac{n-3}{2}}\\
&\le
\frac{1+\frac1{n-1}}2\,t^2
\exp\left[
-\left(\frac12-\frac1{n-1}\right)t
\right].
\end{aligned}
\]
For $n\ge6001$, the function
\[
t\longmapsto
t^2\exp\left[
-\left(\frac12-\frac1{n-1}\right)t
\right]
\]
is decreasing on $[16,\infty)$.

We use the elementary numerical bounds
\begin{equation}\label{eq:exp-constants}
\ee^{\frac{79}{10}}>2600,
\qquad
\ee^6>403,
\qquad
\ee^8>2980,
\qquad
16\ee^{-16}<\frac1{1000}.
\end{equation}
The first three inequalities follow
by taking sufficiently many terms
of the exponential series; summation
through degree $32$ is sufficient.
The last follows from $\ee>2$.

Since
\[
16\left(\frac12-\frac1{n-1}\right)
\ge8-\frac{16}{6000}>\frac{79}{10},
\]
the first inequality in
\eqref{eq:exp-constants} gives
\[
\frac{n\theta\eps^{\frac{n-1}{2}}}{K_0}
<
\frac{128(1+\frac1{6000})}{2600}
<\frac1{20}.
\]

On the other hand, \eqref{eq:C} gives
\[
\frac1u-C_{n-1}(u)
=\frac{\eps^n}{u}.
\]
Using \eqref{eq:K0} and
$t=(n-1)u$, we obtain
\[
\frac{1/u-C_{n-1}(u)}{K_0}
=t\eps^{n-1}\le t\ee^{-t}
<\frac1{1000},
\]
where the final inequality follows
from $t>16$ and
\eqref{eq:exp-constants}.

We have therefore established
\begin{equation}\label{eq:endpoint-small}
\begin{aligned}
\frac{n\theta\eps^{\frac{n-1}{2}}}{K_0}
&<\frac1{20},\\
\frac{1/u-C_{n-1}(u)}{K_0}
&<\frac1{1000}.
\end{aligned}
\end{equation}

Next, since $D=2b-u$ and $0<u/b\le1$,
\[
\frac{u^2D}{b^3}
=\left(\frac ub\right)^2
\left(2-\frac ub\right)\le1.
\]
Indeed, the function
$x^2(2-x)$ is increasing on $[0,1]$
and equals one at $x=1$.

Also, \eqref{eq:Q0-parameters}
gives $m_0\le\eps$. It follows that
\[
\frac1{K_0}
\frac{m_0D}{(n-1)b^3}
=
\frac{m_0}{\eps}
\frac{u^2D}{b^3}\le1.
\]

Finally, using
$b-u=h_0/2=3\eps/((n-1)u)$,
we obtain
\begin{equation}\label{eq:large-b-gap}
\frac{1/u-1/b}{K_0}
=\frac{b-u}{ub}
\frac{(n-1)u^2}{\eps}=\frac3b.
\end{equation}

Combining \eqref{eq:right-moment},
\eqref{eq:endpoint-small}, and
\eqref{eq:large-b-gap}, together with
$m_0\le\eps$ and $u^2D/b^3\le1$,
we conclude that
\begin{equation}\label{eq:large-gap}
\frac{C_{n-1}(u)-R_{n-1}}{K_0}
>
\frac3b-\frac{11}{10}
-\frac1{20}-\frac1{1000}.
\end{equation}

\medskip
\noindent
\emph{Step 2. An upper bound for
$\Err_{n-1}(u)$.}
We now estimate the error in
\eqref{eq:error-definition}.
For any $d\in(0,1)$,
\eqref{eq:error-integral} gives
\begin{equation}\label{eq:large-error-split}
\frac{\Err_{n-1}(u)}{K_0}
\le
\frac1{K_0}\int_0^dO_{n-1}(r)\dd r
+
\frac1{K_0}\int_d^1
\min\{O_{n-1}(r),T_{n-1}(r)\}\dd r.
\end{equation}

By \eqref{eq:error-front} and
\eqref{eq:K0}, the first integral in
\eqref{eq:large-error-split} satisfies
\begin{equation}\label{eq:large-front}
\begin{aligned}
\frac1{K_0}\int_0^dO_{n-1}(r)\dd r
&\le
\frac{40(n-1)^2}
{3(n-2)(n+2)(2-d)^3}\\
&<\frac{40}{3(2-d)^3},
\end{aligned}
\end{equation}
where the strict inequality follows
from $(n-1)^2<(n-2)(n+2)$
for $n\ge4$.

We consider three ranges of $u$.

\smallskip\noindent
\emph{Case 1: $0<u\le1/2$.}
Choose $d=1/2$. Substituting this
value into \eqref{eq:large-front}
gives
\[
\frac1{K_0}\int_0^{1/2}
O_{n-1}(r)\dd r<\frac{320}{81}.
\]

For the integral over $[1/2,1]$,
we use the second bound in
\eqref{eq:error-tail}. Since
\[
\mathcal B_u\left(\frac12\right)
=1-\frac{3u}{4},
\qquad \eps\ge\frac12,
\]
equations \eqref{eq:error-tail} and
\eqref{eq:K0} give
\[
\begin{aligned}
\frac1{K_0}
\int_{1/2}^1
\min\{O_{n-1}(r),T_{n-1}(r)\}\dd r
&\le
\frac{n(n-1)u^2}{2\eps}
\left[
\left(1-\frac{3u}{4}\right)^{\frac{n-1}{2}}
+\left(1-\frac u2\right)^{n-1}
\right]\\
&\le
\left(1+\frac1{n-1}\right)t^2
\left(
\ee^{-\frac{3t}{8}}
+\ee^{-\frac t2}
\right),
\end{aligned}
\]
where the last inequality uses
$\eps\ge1/2$, $t=(n-1)u$,
and $1-x\le\ee^{-x}$.

Both functions
$t^2\ee^{-\frac{3t}{8}}$ and
$t^2\ee^{-\frac t2}$ are decreasing
for $t\ge16$. Since $n\ge6001$,
\eqref{eq:large-u} and
\eqref{eq:exp-constants} imply
\[
\begin{aligned}
\left(1+\frac1{n-1}\right)t^2
\left(
\ee^{-\frac{3t}{8}}+\ee^{-\frac t2}
\right)
&\le
\left(1+\frac1{6000}\right)16^2
\left(\ee^{-6}+\ee^{-8}\right)\\
&<
\left(1+\frac1{6000}\right)256
\left(\frac1{403}+\frac1{2980}\right)
<\frac{73}{100}.
\end{aligned}
\]
Combining this bound with
\eqref{eq:large-error-split} and
\eqref{eq:large-front}, we obtain
\begin{equation}\label{eq:small-u-error}
\frac{\Err_{n-1}(u)}{K_0}
<
\frac{320}{81}+\frac{73}{100},
\qquad 0<u\le\frac12.
\end{equation}

\smallskip
\noindent
\emph{Case 2: $1/2\le u\le4/5$.}

Choose $d=1/4$. Equation
\eqref{eq:large-front} gives
\[
\frac1{K_0}\int_0^{1/4}
O_{n-1}(r)\dd r
<
\frac{40}{3(7/4)^3}
=\frac{2560}{1029}.
\]

\smallskip
\noindent
\emph{Case 3: $4/5\le u<1$.}

Choose $d=1/16$. By
\eqref{eq:large-front},
\[
\frac1{K_0}\int_0^{1/16}
O_{n-1}(r)\dd r
<
\frac{40}{3(31/16)^3}
=\frac{163840}{89373}.
\]

For Cases 2 and 3, the tail integral
admits the same upper bound.
Indeed, for $d=1/4$ and $u\ge1/2$,
\[
\frac{ud(2-d)}2
=\frac{7u}{32}
\ge\frac7{64}>\frac{31}{640}.
\]
For $d=1/16$ and $u\ge4/5$,
\[
\frac{ud(2-d)}2
=\frac{31u}{512}
\ge\frac{31}{640}.
\]
Thus, in both cases,
\begin{equation}\label{eq:large-decay}
\mathcal B_u(d)^{\frac{n-4}{2}}
=
[1-ud(2-d)]^{\frac{n-4}{2}}
\le
\exp\left(-\frac{31(n-4)}{640}\right).
\end{equation}

Using the first bound in
\eqref{eq:error-tail}, followed by
\eqref{eq:K0} and
\eqref{eq:large-decay}, we obtain
\[
\begin{aligned}
\frac1{K_0}\int_d^1
\min\{O_{n-1}(r),T_{n-1}(r)\}\dd r
&\le
\frac56n(n-1)^2u^3(1-d)
\mathcal B_u(d)^{\frac{n-4}{2}}\\
&\le
n(n-1)^2
\exp\left(-\frac{31(n-4)}{640}\right).
\end{aligned}
\]

The function
\[
x\longmapsto
x(x-1)^2
\exp\left(-\frac{31(x-4)}{640}\right)
\]
is decreasing for $x\ge6001$, since
\[
\frac1x+\frac2{x-1}-\frac{31}{640}<0
\qquad(x\ge6001).
\]
At $x=6001$, the exponent satisfies
\[
\frac{31(6001-4)}{640}>280,
\]
while
$
1000\cdot6001\cdot6000^2
<2^{280}<\ee^{280}.
$
Consequently,
\begin{equation}\label{eq:tiny-tail}
n(n-1)^2
\exp\left(-\frac{31(n-4)}{640}\right)
<\frac1{1000},
\qquad n\ge6001.
\end{equation}

Combining \eqref{eq:large-error-split},
\eqref{eq:large-front}, and
\eqref{eq:tiny-tail}, we conclude that
\begin{equation}\label{eq:large-error-cases}
\frac{\Err_{n-1}(u)}{K_0}
<
\begin{cases}
\displaystyle
\frac{2560}{1029}+\frac1{1000},
&\displaystyle \frac12\le u\le\frac45,\\[8pt]
\displaystyle
\frac{163840}{89373}+\frac1{1000},
&\displaystyle \frac45\le u<1.
\end{cases}
\end{equation}

\medskip
\noindent
\emph{Step 3. Comparison of the bounds.}

By the definition of $b$ in
\eqref{eq:Q0-parameters} and the choice
$\kappa=6/(n-1)$, we have
\[
b=b(u)
=u-\frac3{n-1}+\frac3{(n-1)u}.
\]
Differentiating twice gives
\[
b''(u)=\frac6{(n-1)u^3}>0.
\]
Thus $b(u)$ is convex on $(0,1)$.

By \eqref{eq:large-u}, we have
$u>16/(n-1)$. On
$[16/(n-1),1/2]$, convexity gives
\[
\begin{aligned}
b(u)
&\le\max\left\{
b\left(\frac{16}{n-1}\right),
b\left(\frac12\right)
\right\}\\
&=\max\left\{
\frac3{16}+\frac{13}{n-1},
\frac12+\frac3{n-1}
\right\}\\
&<\frac{501}{1000},
\end{aligned}
\]
where the last inequality uses
$n-1\ge6000$.

Similarly, on $[1/2,4/5]$,
\[
\begin{aligned}
b(u)
&\le\max\left\{
b\left(\frac12\right),
b\left(\frac45\right)
\right\}\\
&=\max\left\{
\frac12+\frac3{n-1},
\frac45+\frac3{4(n-1)}
\right\}\\
&<\frac{801}{1000}.
\end{aligned}
\]

Finally, for $4/5\le u<1$,
the inequality $(n-1)u>16$ gives
\[
b=u+\frac{3(1-u)}{(n-1)u}
\le u+(1-u)=1.
\]

We have therefore established
\begin{equation}\label{eq:large-b-bounds}
b<
\begin{cases}
\displaystyle \frac{501}{1000},
&0<u\le\frac12,\\[5pt]
\displaystyle \frac{801}{1000},
&\frac12\le u\le\frac45,
\end{cases}
\qquad
b\le1\quad\text{if }\frac45\le u<1.
\end{equation}

Since $K_0>0$, equation
\eqref{eq:large-gap} gives
\begin{equation}\label{eq:large-final-gap}
\frac{
C_{n-1}(u)-\Err_{n-1}(u)-R_{n-1}
}{K_0}>
\frac3b-\frac{11}{10}
-\frac1{20}-\frac1{1000}
-\frac{\Err_{n-1}(u)}{K_0}.
\end{equation}

Substituting \eqref{eq:small-u-error},
\eqref{eq:large-error-cases}, and
\eqref{eq:large-b-bounds} into
\eqref{eq:large-final-gap}, we obtain
the explicit positive lower bounds
in Table~\ref{tab:large}.
In particular,
\[
C_{n-1}(u)-\Err_{n-1}(u)>R_{n-1}.
\]
This contradicts
\eqref{eq:quadratic-comparison},
which follows from
\eqref{eq:second} and the choice of
$h_0<h$. Hence \eqref{eq:polar}
and \eqref{eq:second} cannot hold
simultaneously for $n\ge6001$.
\end{proof}

\begin{table}[htbp]
\centering
\small
\caption{Lower bounds for
$(C_{n-1}(u)-\Err_{n-1}(u)-R_{n-1})/K_0$
when $n\ge6001$.}
\label{tab:large}
\renewcommand{\arraystretch}{1.25}
\begin{tabular}{@{}cccc@{}}
\toprule
Range of $u$
& \shortstack{Upper bound\\for $b$}
& \shortstack{Upper bound\\for $\Err_{n-1}(u)/K_0$}
& Lower bound\\
\midrule
$(0,1/2]$
& $501/1000$
& $320/81+73/100$
& $2115713/13527000$\\
$[1/2,4/5]$
& $801/1000$
& $2560/1029+1/1000$
& $1207336/11447625$\\
$[4/5,1)$
& $1$
& $163840/89373+1/1000$
& $165163/11171625$\\
\bottomrule
\end{tabular}
\end{table}

\subsection{The remaining degrees}

\begin{lem}[Exclusion of the remaining degrees]
\label{lem:finite}
For every integer $4\le n\le6000$,
there are no $0<u<1$ and $0<h<1$
satisfying both \eqref{eq:polar}
and \eqref{eq:second}.
\end{lem}

\begin{proof}
Suppose that $4\le n\le6000$ and
that $0<u<1$ and $0<h<1$ satisfy
\eqref{eq:polar} and \eqref{eq:second}.

By Lemma~\ref{lem:polar},
\begin{equation}\label{eq:finite-threshold}
u\ge\tau_{n-1}
=\left(\frac{n+1}{2}\right)^{\frac2{n-1}}-1.
\end{equation}

The interval certificates in the accompanying repository~\cite{Zha26e},
described in Appendix~\ref{app:verification}, prove
the following alternative for every
integer $4\le n\le6000$ and every
$u\in[\tau_{n-1},1)$:
either the hypotheses of
Lemma~\ref{lem:local} are satisfied,
or a certified value of $\kappa>0$
satisfies
\begin{equation}\label{eq:finite-certificate}
h_0=\frac{\kappa\eps}{u}<h,
\qquad
\frac{\kappa}{u}\le4,
\qquad
C_{n-1}(u)-\Err_{n-1}(u)>R_{n-1}.
\end{equation}

If Lemma~\ref{lem:local} applies,
then \eqref{eq:polar} and
\eqref{eq:second} are incompatible,
contrary to our assumption.

Otherwise, the first two conditions
in \eqref{eq:finite-certificate}
ensure that the comparison
\eqref{eq:quadratic-comparison}
is valid. Thus
\[
C_{n-1}(u)-\Err_{n-1}(u)
\le R_{n-1},
\]
contradicting the third inequality
in \eqref{eq:finite-certificate}.

For completeness, the certificates
cover the parameter interval $[0,1]$
by closed intervals. Intervals below
the threshold in
\eqref{eq:finite-threshold} are excluded
by Lemma~\ref{lem:polar}, while the
endpoint $u=1$ is excluded by the
argument following that lemma.
The finite verification consists of
$261$ consecutive ranges of degrees
and $20107$ terminal intervals,
with the counts recorded in
Table~\ref{tab:certificates}.
Every accepted interval comparison
uses the outward-rounded integer
arithmetic specified in
Appendix~\ref{app:arithmetic}.

It follows that \eqref{eq:polar}
and \eqref{eq:second} cannot hold
simultaneously for $4\le n\le6000$.
\end{proof}

\Needspace{8\baselineskip}

We can now complete the proof of
Conjecture~\ref{conj:Schmeisser}.

\begin{proof}[Proof of Conjecture~\ref{conj:Schmeisser}]
Let $p$ be a polynomial of degree
$n\ge2$ whose zeros lie in
$\overline{\mathbb D}$.

If $n=2$, the unique critical point
is the midpoint of the two zeros.
Every point in their convex hull
is therefore at distance at most
half the distance between the zeros
from the critical point, and this
distance is at most one.

Suppose now that $n\ge3$.
Let $\xi\in H(p)$.
By Lemma~\ref{lem:diameter-cover},
there is a zero $a$ of $p$ such that
$\xi\in D_a$.

If $a$ is a multiple zero, then
$p'(a)=0$. Since $\xi\in D_a$,
\[
|\xi-a|
\le\left|\xi-\frac a2\right|
+\frac{|a|}{2}
\le|a|\le1.
\]
If $a=0$, then $D_a=\{0\}$,
so $\xi=0$. The Gauss--Lucas theorem
implies that every critical point
of $p$ lies in $\overline{\mathbb D}$.
Hence $\xi$ is again within distance
one of a critical point.

It remains to consider the case where
$a$ is a nonzero simple zero.
After rotating the plane, we may
assume $a\in(0,1]$.

If $a=1$, Lemma~\ref{lem:boundary-zero}
shows that \eqref{eq:matrix} holds.
If $0<a<1$, suppose that
\eqref{eq:matrix} fails.
Lemma~\ref{lem:reduction} then gives
$0<u\le1$ and $0<h<1$ satisfying
\eqref{eq:polar} and \eqref{eq:second}.

The case $u=1$ is impossible:
by \eqref{eq:C} and
\eqref{eq:error-definition},
the left-hand side of
\eqref{eq:second} equals one,
whereas its right-hand side is
strictly less than
\[
n\int_0^1(1-t)^{n-1}\dd t=1
\]
because $h>0$.

If $n=3$, the first inequality in
\eqref{eq:polar-bounds} gives
$u\ge\tau_2=1$, which is impossible
by the exclusion of $u=1$.
For $n\ge4$, Lemma~\ref{lem:large}
excludes all degrees $n\ge6001$,
while Lemma~\ref{lem:finite} excludes
$4\le n\le6000$.

Thus \eqref{eq:matrix} holds at every
nonzero simple zero. Applying
Lemma~\ref{lem:geometry}, we conclude
that every point of $D_a$, and in
particular $\xi$, lies within distance
one of a critical point of $p$.

Finally, the constant one is sharp.
For $p(z)=z^n-1$, the only critical
point is the origin, while every
zero of $p$ lies on the unit circle.
Hence the distance from any zero
of $p$ to its nearest critical
point is exactly one.
\end{proof}

\appendix
\section{The finite verification}\label{app:verification}

We describe the interval inequalities used to exclude degrees
$4\le n\le6000$, together with the arithmetic rules of the
independent verifier. The functions $C_{n-1}$, $\Err_{n-1}$,
$Q_0$, and $R_m$ are defined in \eqref{eq:C},
\eqref{eq:error-definition}, \eqref{eq:Q0-definition}, and
\eqref{eq:quadratic}, respectively. Each certificate establishes
strict inequalities on entire closed intervals of $u$; verification
at isolated sample points would not suffice.

\subsection{An endpoint estimate for quartics}

\begin{lem}[The quartic endpoint]\label{lem:local}
For $n=4$ and $0<1-u\le1/100$, the inequalities
\eqref{eq:polar} and \eqref{eq:second} cannot hold simultaneously.
\end{lem}

\begin{proof}
Suppose that $0<\eps:=1-u\le1/100$ and that
\eqref{eq:polar} holds. In particular, $u\ge99/100$.
Choose $\kappa=3$ and set $h_0=3\eps/u$. Then
\[
0<h_0\le\frac1{33}<\frac1{20},
\qquad \frac\kappa u\le\frac{100}{33}<4.
\]
We first prove that $h_0<h$.

For $0\le t\le1$, put $X(t)=(1+u)t-1$.
Since $-1\le X(t)\le u$, Taylor's theorem applied to
$x\mapsto(1+x)^{\frac32}$ on $[-h_0,h_0u]$ gives
\[
(1+h_0X(t))^{\frac32}
\le 1+\frac32h_0X(t)
+\frac{3h_0^2X(t)^2}{8\sqrt{1-h_0}}.
\]
As
\[
\int_0^1X(t)\dd t=-\frac\eps2,
\qquad
\int_0^1X(t)^2\dd t
=\frac{u^2-u+1}{3}\le\frac13,
\]
the definition of $P_3$ yields
\begin{equation}\label{eq:quartic-polar-test}
P_3(u,h_0)-1
\le-\frac34h_0\eps
+\frac{h_0^2}{8\sqrt{1-h_0}}<0.
\end{equation}
To justify the strict inequality, note that it is equivalent to
$h_0<6\eps\sqrt{1-h_0}$, or
$2u\sqrt{1-h_0}>1$. The inequality
$2u\sqrt{1-h_0}>1$ follows from
$u\ge99/100$ and $h_0\le1/33$.
If $h\le h_0$, convexity of $P_3(u,\cdot)$, together
with $P_3(u,0)=1$, would imply
\[
P_3(u,h)\le
\left(1-\frac h{h_0}\right)P_3(u,0)
+\frac h{h_0}P_3(u,h_0)<1,
\]
contrary to \eqref{eq:polar} and \eqref{eq:quartic-polar-test}.
Therefore $h_0<h$.

The restriction $\kappa/u<4$ ensures $m_0\ge0$ in
\eqref{eq:Q0-parameters}. By \eqref{eq:Q0-square},
$\sqrt{Q_0(t)}=\sqrt{m_0+D(t-b/D)^2}$ is a convex
function of $t$. Since $Q_0(0)=1$ and $Q_0(1)=\eps$,
convexity gives
\[
\sqrt{Q_0(t)}\le1-t+t\sqrt\eps,
\qquad 0\le t\le1.
\]
Integrating and using \eqref{eq:quadratic} with $m=1$, we obtain
\begin{equation}\label{eq:quartic-R1}
R_1=2\int_0^1\sqrt{Q_0(t)}\dd t\le1+\sqrt\eps.
\end{equation}
Moreover, \eqref{eq:theta-definition} and
$h_0=3\eps/u$ give
\begin{equation}\label{eq:quartic-theta}
\frac\theta\eps
=\frac{3}{2(u^2+3\eps)}
=\frac{3}{2(1+\eps+\eps^2)}
>\frac{37}{25}.
\end{equation}
Here $\eps\le1/100$ is used in the final comparison.
Since $m_0\le\eps$ by \eqref{eq:Q0-parameters},
\eqref{eq:recurrence} with $m=3$, together with
\eqref{eq:quartic-R1} and \eqref{eq:quartic-theta}, yields
\begin{equation}\label{eq:quartic-R3}
\begin{aligned}
R_3
&=1-\theta(1-\eps^{\frac32})+\frac32m_0R_1\\
&\le1-\theta(1-\eps^{\frac32})
+\frac32\eps(1+\sqrt\eps)\\
&<1+\eps\left[
-\frac{37}{25}\left(1-\frac1{1000}\right)
+\frac32\left(1+\frac1{10}\right)\right]\\
&<1+\frac7{40}\eps.
\end{aligned}
\end{equation}

For $n=4$, the coefficients in \eqref{eq:error-finite}
are $a_{3,2}=3/2$ and $a_{3,3}=1$. Hence
\begin{equation}\label{eq:quartic-error}
\begin{aligned}
\frac{\Err_3(u)}\eps
&\le\frac u2(1+3\eps)+u^{\frac32}\sqrt\eps\\
&\le\frac12\left(1+\frac3{100}\right)
+\frac1{10}=\frac{123}{200}.
\end{aligned}
\end{equation}
Finally, \eqref{eq:C} gives
$C_3(u)=1+\eps+\eps^2+\eps^3\ge1+\eps$.
Combining \eqref{eq:quartic-R3} and
\eqref{eq:quartic-error}, we find that
\[
C_3(u)-\Err_3(u)-R_3
>\left(1-\frac{123}{200}-\frac7{40}\right)\eps
=\frac{21}{100}\eps>0.
\]
Because $h_0<h$ and $\kappa/u<4$, the comparison
\eqref{eq:quadratic-comparison} applies and gives the opposite
inequality $C_3(u)-\Err_3(u)\le R_3$.
This contradiction proves the lemma.
\end{proof}

\subsection{A normalized recurrence}

For $0<u<1$, define
\begin{equation}\label{eq:J-definition}
\mathcal J_m:=\frac{1/u-R_m}{\eps},
\qquad m\ge0,
\end{equation}
and retain the normalization\footnote{This is the horizontal
distance $\theta$ from the
vertex to the right endpoint, normalized by $\eps=Q_0(1)$. The symbol
$t_0$ denotes this ratio, not the location of the vertex.}
\begin{equation}\label{eq:t0-definition}
t_0:=\frac\theta\eps=\frac\kappa{2uD}.
\end{equation}
The identities $h_0=\kappa\eps/u$ and
\eqref{eq:Q0-parameters} also give
\begin{equation}\label{eq:m0-normalization}
\frac{m_0}\eps
=1-\frac{\kappa^2\eps}{4u^2D},
\qquad 0\le\frac{m_0}\eps\le1.
\end{equation}

Substituting
$R_{m-2}=1/u-\eps\mathcal J_{m-2}$
into \eqref{eq:recurrence}, then using
\eqref{eq:t0-definition} and \eqref{eq:m0-normalization},
we obtain, for every integer $m\ge2$,
\begin{equation}\label{eq:J-recurrence}
\mathcal J_m
=-\frac1{(m-1)u}
+\frac{m}{m-1}\frac{\kappa^2\eps}{4u^3D}+t_0(1-\eps^{\frac m2})
+\frac{m}{m-1}m_0\mathcal J_{m-2}.
\end{equation}
Because $m_0\ge0$, the coefficient of
$\mathcal J_{m-2}$ in \eqref{eq:J-recurrence}
is nonnegative; interval lower bounds can therefore be
propagated through the recurrence.

By \eqref{eq:F-comparison} and the three inequalities
in \eqref{eq:F-elementary},
\begin{equation}\label{eq:J-tail}
\begin{aligned}
\mathcal J_m&\ge\frac1u-m,\\
\mathcal J_m&\ge-\frac1{u\eps},\\
\mathcal J_m&\ge\frac1u-3
\qquad\text{if }\eps\le\frac14,
\end{aligned}
\qquad m\ge2.
\end{equation}
For $m=2$, direct integration of
$Q_0(t)=1-2bt+Dt^2$ in \eqref{eq:quadratic}
gives $R_2=1+2\eps-h_0/2$, and hence
\begin{equation}\label{eq:J2}
\mathcal J_2=\frac1u-2+\frac\kappa{2u}.
\end{equation}

For $m=3$, a sharper estimate comes from
an explicit upper bound on $R_1$.
Write
\[
\varphi(t):=\sqrt{Q_0(t)}-(1-t),
\qquad 0\le t\le1.
\]
Since $Q_0(t)=m_0+D(t-b/D)^2$ and $m_0\ge0$,
$\varphi$ is convex.\footnote{The function $\varphi$ is the signed
vertical gap between $y=\sqrt{Q_0(t)}$ and $y=1-t$. Its integral is their
signed area, estimated by trapezoids in \eqref{eq:R1-trapezoid}.}
The identity
\[
Q_0(t)-(1-t)^2
=\eps\left[t(2-t)-\frac\kappa u t(1-t)\right]
\]
allows us to evaluate $\varphi$ without subtracting
nearly equal square roots:
\begin{equation}\label{eq:trapezoid}
\begin{aligned}
r_j&:=\frac j{16}\qquad(0\le j\le16),\\
\varphi(t)&=
\frac{\eps\left[t(2-t)-(\kappa/u)t(1-t)\right]}
{\sqrt{Q_0(t)}+1-t},\\
V_1&:=\frac{\varphi(0)+\varphi(1)
+2\sum_{j=1}^{15}\varphi(r_j)}{16}.
\end{aligned}
\end{equation}
Here $\varphi(0)=0$ and $\varphi(1)=\sqrt\eps$.
For a convex function, the trapezoidal rule on a
uniform partition gives an upper bound for its integral.
Applying that rule to the sixteen panels
$[r_j,r_{j+1}]$ and using \eqref{eq:trapezoid},
we obtain
\begin{equation}\label{eq:R1-trapezoid}
R_1-1=2\int_0^1\varphi(t)\dd t\le V_1.
\end{equation}
The sign of $\varphi$ is immaterial to the
trapezoidal inequality.

From \eqref{eq:recurrence} with $m=3$,
\eqref{eq:J-definition}, \eqref{eq:t0-definition},
\eqref{eq:m0-normalization}, and \eqref{eq:R1-trapezoid},
we obtain
\begin{equation}\label{eq:J3}
\mathcal J_3\ge\frac1u+t_0(1-\eps^{\frac32})-\frac32
\left(1-\frac{\kappa^2\eps}{4u^2D}\right)(1+V_1).
\end{equation}
The substitution of $R_1\le1+V_1$ is valid
because $m_0/\eps\ge0$ by
\eqref{eq:m0-normalization}.

When $(n-1)\alpha\ge90$, combining
\eqref{eq:right-moment}, \eqref{eq:J-definition},
\eqref{eq:t0-definition}, and
\eqref{eq:m0-normalization} gives
\begin{equation}\label{eq:J-moment}
\mathcal J_{n-1}\ge\frac\kappa{2u^2b}
-\frac{11}{10}\frac{D}{(n-1)b^3}
\left(1-\frac{\kappa^2\eps}{4u^2D}\right)
-nt_0\eps^{\frac{n-1}{2}}.
\end{equation}
Indeed, the first term follows from
$(1/u-1/b)/\eps=\kappa/(2u^2b)$.

The functions \href{https://github.com/zhangteng2000/schmeisser-conjecture/blob/9e958aec9d41c387d19a2485fa771895f09eda60/work/verify_mass_scalar_integer_v4.py#L85-L101}{\nolinkurl{Checker.recurrence}},
\href{https://github.com/zhangteng2000/schmeisser-conjecture/blob/9e958aec9d41c387d19a2485fa771895f09eda60/work/verify_mass_scalar_integer_v4.py#L72-L84}{\nolinkurl{Checker.base_three}}, and
\href{https://github.com/zhangteng2000/schmeisser-conjecture/blob/9e958aec9d41c387d19a2485fa771895f09eda60/work/verify_mass_scalar_integer_v4.py#L68-L71}{\nolinkurl{Checker.moment}}
in \nolinkurl{work/verify_mass_scalar_integer_v4.py}
implement \eqref{eq:J-recurrence}, \eqref{eq:J3}, and
\eqref{eq:J-moment}, respectively. The terminal comparison is
performed by \href{https://github.com/zhangteng2000/schmeisser-conjecture/blob/9e958aec9d41c387d19a2485fa771895f09eda60/work/verify_mass_scalar_integer_v4.py#L110-L127}{\nolinkurl{Checker.check}}.

The strict inequality needed in the finite verification
is $C_{n-1}(u)-\Err_{n-1}(u)>R_{n-1}$.
By \eqref{eq:C} and \eqref{eq:J-definition},
that inequality is equivalent to
\begin{equation}\label{eq:acceptance}
\mathcal J_{n-1}>
\frac{\Err_{n-1}(u)}\eps
+\frac{\eps^{n-1}}u.
\end{equation}

\Needspace{6\baselineskip}
\subsection{Bounds over a range of degrees}\label{app:degree-ranges}

Fix integers $4\le n_-\le n_+\le6000$.
The verifier processes the entire range
$n_-\le n\le n_+$ in one certificate.
By \eqref{eq:polar-bounds}, any admissible $u$ satisfies
$u\ge\tau_{n-1}$. To find a lower threshold valid
throughout the degree range, observe that
\[
\tau_{n-1}=
\exp\left(
\frac{\log(1+(n-1)/2)}{(n-1)/2}
\right)-1.
\]
Since $x\mapsto\log(1+x)/x$ decreases on $(0,\infty)$,
$\tau_{n-1}$ decreases with $n$. Thus
\begin{equation}\label{eq:degree-threshold}
u\ge\tau_{n_+-1}
\qquad(n_-\le n\le n_+).
\end{equation}

The default choice is $\kappa=6/(n_+-1)$.
For $n\le n_+$, \eqref{eq:polar-bounds} gives
\[
\frac{\kappa\eps}{u}
=\frac{6\eps}{(n_+-1)u}
\le\frac{6\eps}{(n-1)u}<h.
\]
The lower bound in \eqref{eq:degree-threshold}
and the estimate $(n_+-1)\tau_{n_+-1}>3/2$
ensure $\kappa/u<4$.
The verifier may also use
\[
\kappa=\frac{j}{4096(n_+-1)},
\qquad j\in\mathbb Z,\quad j\ge24576.
\]
For every alternative choice, the verifier certifies
$0<h_0\le1$, $\kappa/u\le4$, and
$P_{n-1}(u,h_0)<1$, the last condition being
computed from \eqref{eq:exact-polar}.
Convexity of $P_{n-1}(u,\cdot)$ and
\eqref{eq:polar} then imply $h_0<h$.

For fixed $u$ and $h_0\in(0,1]$, the integrand
$[1-h_0+h_0(1+u)t]^{\frac{n-1}{2}}$
is convex in the real variable $n$ whenever its base
is positive. It follows that
\begin{equation}\label{eq:degree-polar-chord}
P_{n-1}(u,h_0)
\le\max\{P_{n_--1}(u,h_0),P_{n_+-1}(u,h_0)\}
\quad(n_-\le n\le n_+).
\end{equation}
Thus checking the two endpoint degrees in
\eqref{eq:degree-polar-chord} certifies the
inequality $P_{n-1}(u,h_0)<1$ for every
integer degree in the range.

The verifier uses \eqref{eq:error-finite} only for
a single degree $n_-=n_+\le129$.
For the geometric estimate \eqref{eq:error-geometric},
it takes $k_*=\min\{10,n_--1\}$ and replaces
$n$ by $n_-$ in the coefficients
$((k-1)/(n-2))^{\frac k2}$.
This replacement increases the coefficients for
$n\ge n_-$. If $n_+-1>k_*$, the geometric
remainder in \eqref{eq:error-geometric} is retained;
that estimate is used only when $\eps/u<1$.

The coefficient of $1/u^2$ in
\eqref{eq:error-whole} and \eqref{eq:error-front}
is proportional to
$
{n-1}/{[(n-2)(n+2)]}.
$
Its derivative as a function of real $n>2$ is
\[
-\frac{n^2-2n+4}{(n^2-4)^2}<0.
\]
Hence the coefficient is bounded above throughout
the degree range by its value at $n=n_-$.
The split points used with \eqref{eq:error-tail} are
$d\in\{1/2,1/4,1/16,3/4\}$.

To control expressions containing a polynomial factor
and an exponential factor in $n$, we use the
following convexity estimate.

\begin{lem}[A product bound for convex functions]
\label{lem:chord}
Let $k$ and $f$ be convex functions on $[A,B]$,
where $A<B$, with $k>0$ and $f\ge0$.
Write $k_-=k(A)$, $k_+=k(B)$,
$f_-=f(A)$, and $f_+=f(B)$.
Then, for every $x\in[A,B]$,
\begin{equation}\label{eq:chord}
k(x)f(x)
\le\max\{k_-f_-,k_+f_+\}
\left(1+\frac{(k_+-k_-)^2}{4k_-k_+}\right).
\end{equation}
\end{lem}

\begin{proof}
Write $x=(1-\lambda)A+\lambda B$ with
$0\le\lambda\le1$, and put
$M=\max\{k_-f_-,k_+f_+\}$.
Convexity gives
\[
k(x)\le(1-\lambda)k_-+\lambda k_+,
\qquad
f(x)\le M\left(\frac{1-\lambda}{k_-}
+\frac\lambda{k_+}\right).
\]
Both factors are nonnegative, and their product
is at most
\[
M\left[
1+\frac{\lambda(1-\lambda)(k_+-k_-)^2}
{k_-k_+}
\right].
\]
Since $\lambda(1-\lambda)\le1/4$,
\eqref{eq:chord} follows.
\end{proof}

For the $O_{n-1}$ term in \eqref{eq:OT},
Lemma~\ref{lem:chord} is applied with
$k(n)=n(n-1)$. For $T_{n-1}$, the factor
is $k(n)=n$. After factoring out the positive
quantities independent of $n$, the remaining
functions are sums or products of exponentials
in $n$, hence nonnegative and convex.

In the rectangle-rule estimate of
\eqref{eq:error-integral}, partition $[0,1]$
into panels $[t_-,t_+]$. Since
$\mathcal B_u(t)$ and $1-ut$ decrease on $[0,1]$,
\eqref{eq:OT} gives, whenever $1-ut_+>0$,
\begin{equation}\label{eq:panel-bounds}
\begin{aligned}
O_{n-1}(t)
&\le\frac56n(n-1)u\eps
\frac{t_+^2}{1-ut_+}
\mathcal B_u(t_-)^{\frac{n-2}{2}},\\
T_{n-1}(t)
&\le n\left[
\mathcal B_u(t_-)^{\frac{n-1}{2}}
+(1-ut_-)^{n-1}\right],
\quad t_-\le t\le t_+.
\end{aligned}
\end{equation}
For each panel, multiply the smaller certified upper
bound from \eqref{eq:panel-bounds} by its width
$t_+-t_-$. If $1-ut_+$ has no certified positive
lower bound, the verifier uses only the bound
for $T_{n-1}$. The resulting panel bounds sum to
an upper bound for the integral in
\eqref{eq:error-integral}, because the panels
cover $[0,1]$ without gaps.

For the normalized recurrence \eqref{eq:J-recurrence},
after $\ell$ steps the index $m$ lies in
\[
[n_--1-2\ell,n_+-1-2\ell].
\]
At each step, the verifier encloses the entire
interval of possible values of $m/(m-1)$.
A lower bound for $\mathcal J_{m-2}$ may be negative,
so replacing $m/(m-1)$ by a single endpoint
without interval multiplication would not be valid.
For $0\le\eps\le1$, the functions
$m\mapsto\eps^{\frac m2}$ are bounded by the
values at the endpoint exponents, with the usual
continuous interpretation at $\eps=0$.
The general starting bounds are \eqref{eq:J-tail};
\eqref{eq:J2} and \eqref{eq:J3} are used only when
$n_-=n_+$.

The functions \href{https://github.com/zhangteng2000/schmeisser-conjecture/blob/9e958aec9d41c387d19a2485fa771895f09eda60/work/verify_mass_scalar_integer_v4.py#L102-L109}{\nolinkurl{Checker.polar_candidate}},
\href{https://github.com/zhangteng2000/schmeisser-conjecture/blob/9e958aec9d41c387d19a2485fa771895f09eda60/work/verify_mass_scalar_integer_v4.py#L33-L54}{\nolinkurl{Checker.error}}, and
\href{https://github.com/zhangteng2000/schmeisser-conjecture/blob/9e958aec9d41c387d19a2485fa771895f09eda60/work/verify_mass_scalar_integer_v4.py#L55-L67}{\nolinkurl{Checker.riemann}}
implement the tests for $\kappa$, the error bounds, and the panel
estimates above. In the repository, the field
\texttt{N\_range} specifies the range of critical-point counts $n-1$.

The bound \eqref{eq:J-moment} is used only after
checking $(n_--1)\alpha\ge90$. The negative term
containing $1/(n-1)$ is bounded using $n=n_-$,
while the final negative term is bounded in magnitude by
$n_+t_0\eps^{\frac{n_--1}{2}}$.
Finally, because $0\le\eps\le1$, the remainder
$\eps^{n-1}/u$ in \eqref{eq:acceptance} is at most
$\eps^{n_--1}/u$ throughout the degree range.

\subsection{Integer arithmetic and coverage}\label{app:arithmetic}

Each certificate specifies a degree range
$n_-\le n\le n_+$, equivalently a critical-point
count range $n_--1\le n-1\le n_+-1$, together
with root intervals in $u$ and a binary subdivision
tree for every root interval. Every subdivision point
is a dyadic rational strictly inside the parent interval.
The two closed children share the subdivision point
and cover the parent exactly; no endpoint is discarded.

The verifier accepts a terminal interval only if
one of the following independently checkable rules applies:
\begin{enumerate}
\item The upper endpoint of the interval is smaller than
an outward-rounded lower bound for $\tau_{n_+-1}$,
so \eqref{eq:degree-threshold} rules out the entire interval.
\item The degree range consists only of $n=4$ and
$1-u\le1/100$ throughout the interval. For $u<1$,
Lemma~\ref{lem:local} rules out simultaneous satisfaction
of \eqref{eq:polar} and \eqref{eq:second}.
\item A recorded value of $\kappa$ satisfies the
admissibility tests of Appendix~\ref{app:degree-ranges},
and the certified lower bound for $\mathcal J_{n-1}$
is strictly greater than the certified upper bound
for the right-hand side of \eqref{eq:acceptance}.
\end{enumerate}

The verifier may first intersect a candidate interval
with the admissible region $u\ge\tau_{n_+-1}$;
\eqref{eq:degree-threshold} already excludes the
part below that threshold. The endpoint $u=1$ is
excluded analytically by the argument following
Lemma~\ref{lem:polar}. Expressions containing
$1/\eps$ are evaluated only for $u<1$.
On intervals meeting $u=1$, the verifier uses
bounds that remain finite at $\eps=0$ or omits
the inapplicable bound.

The subdivision trees were generated using
$112$-bit Arb interval arithmetic
\cite{Joh17}. The independent verifier neither
imports Arb nor calls the generator. Its intervals
have the fixed-point form
\begin{equation}\label{eq:fixed-point-interval}
\left[\frac{\ell_*}{s_*},\frac{r_*}{s_*}\right],
\qquad s_*:=2^{192},
\end{equation}
where $\ell_*$ and $r_*$ are integers.
Addition and subtraction of endpoints on the grid
in \eqref{eq:fixed-point-interval} are exact.
Multiplication and division use integer floor and
ceiling to enclose their exact results. Square-root
bounds are obtained from integer square roots,
with outward adjustment of the upper endpoint.
All powers used by the verifier have nonnegative
integer or half-integer exponents.
The minimum of certified upper bounds is still an
upper bound, and the maximum of certified lower
bounds is still a lower bound.
The arithmetic routines may impose the established
constraints
\[
\eps\ge0,\qquad
0\le1-\frac{\kappa^2\eps}{4u^2D}\le1,
\qquad Q_0(t)\ge0,
\]
provided the admissibility conditions on $\kappa$
are satisfied; the middle inequality follows from
\eqref{eq:m0-normalization}.

The class \href{https://github.com/zhangteng2000/schmeisser-conjecture/blob/9e958aec9d41c387d19a2485fa771895f09eda60/work/mass_integer_arithmetic.py#L10-L83}{\nolinkurl{I}}
in \href{https://github.com/zhangteng2000/schmeisser-conjecture/blob/9e958aec9d41c387d19a2485fa771895f09eda60/work/mass_integer_arithmetic.py}{\nolinkurl{work/mass_integer_arithmetic.py}} implements this
interval arithmetic; its function
\href{https://github.com/zhangteng2000/schmeisser-conjecture/blob/9e958aec9d41c387d19a2485fa771895f09eda60/work/mass_integer_arithmetic.py#L86-L98}{\nolinkurl{threshold}}
computes the threshold enclosure described next.

The threshold $\tau_{n-1}$ is enclosed without
floating-point assumptions. For each relevant
integer degree $n$, the verifier finds an integer
$z$ satisfying
\begin{equation}\label{eq:integer-root}
z^{n-1}\le\frac{(n+1)^2s_*^{n-1}}4
<(z+1)^{n-1}.
\end{equation}
Since $\tau_{n-1}=((n+1)/2)^{\frac2{n-1}}-1$,
\eqref{eq:integer-root} implies
\begin{equation}\label{eq:threshold-interval}
\frac z{s_*}-1\le\tau_{n-1}
<\frac{z+1}{s_*}-1.
\end{equation}
The enclosure in \eqref{eq:threshold-interval}
is used in the threshold rejection test of rule~1.
A floating-point approximation is used only to
locate a candidate integer $z$; both comparisons
in \eqref{eq:integer-root} are then checked exactly
with integers. Dyadic subdivision points are parsed
as exact rationals, so their floating-point
approximations do not affect strict interval tests.

The verifier traverses every subdivision tree from
its root to its terminal intervals. An end marker is
accepted only when the stack of unprocessed
intervals is empty and the certificate contains no
additional data. A separate check confirms that the
degree ranges are consecutive and cover precisely
the critical-point counts $3,\ldots,5999$.
The recorded counts appear in
Table~\ref{tab:certificates}.

\begin{table}[htbp]
\centering
\small
\caption{Degree ranges and interval certificates
for the finite verification.}
\label{tab:certificates}
\begin{tabular}{@{}lrrr@{}}
\toprule
Number of critical points $n-1$
& Degree ranges & Tree nodes & Terminal intervals\\
\midrule
$3$--$19$ & 17 & 27319 & 13668\\
$20$--$128$ & 109 & 4977 & 2543\\
$129$--$5999$ & 135 & 7657 & 3896\\
\midrule
Total & 261 & 39953 & 20107\\
\bottomrule
\end{tabular}
\end{table}

The function \href{https://github.com/zhangteng2000/schmeisser-conjecture/blob/9e958aec9d41c387d19a2485fa771895f09eda60/work/verify_mass_scalar_integer_v4.py#L129-L155}{\nolinkurl{verify}}
checks every subdivision tree. The consecutive degree ranges and
the totals in Table~\ref{tab:certificates} are checked by
\href{https://github.com/zhangteng2000/schmeisser-conjecture/blob/9e958aec9d41c387d19a2485fa771895f09eda60/verify_all.py#L29-L91}{\nolinkurl{main}}
in \href{https://github.com/zhangteng2000/schmeisser-conjecture/blob/9e958aec9d41c387d19a2485fa771895f09eda60/verify_all.py}{\nolinkurl{verify_all.py}}, using
\href{https://github.com/zhangteng2000/schmeisser-conjecture/blob/9e958aec9d41c387d19a2485fa771895f09eda60/work/mass_proof_manifest.json}{\nolinkurl{work/mass_proof_manifest.json}}.

The counts in Table~\ref{tab:certificates}
are consistent with binary subdivision forests:
for each row, the number of nodes equals twice
the number of terminal intervals minus the number
of root trees (degree ranges).

\subsection{Reproduction}

The code and certificates are available in the public repository
\cite{Zha26e}. All links in this appendix refer to commit
\href{https://github.com/zhangteng2000/schmeisser-conjecture/tree/9e958aec9d41c387d19a2485fa771895f09eda60}{\texttt{9e958aec9d41}},
which fixes the version used for the verification.
The repository contains the certificate files, their generator,
the independent integer verifier, the arithmetic kernel,
\href{https://github.com/zhangteng2000/schmeisser-conjecture/blob/9e958aec9d41c387d19a2485fa771895f09eda60/work/mass_proof_manifest.json}{\nolinkurl{work/mass_proof_manifest.json}},
and \href{https://github.com/zhangteng2000/schmeisser-conjecture/blob/9e958aec9d41c387d19a2485fa771895f09eda60/SHA256SUMS.txt}{\nolinkurl{SHA256SUMS.txt}}.
Running \texttt{python verify\_all.py} checks the
recorded file hashes, replays all $261$ degree ranges,
and verifies the $33$ rational comparisons in
\href{https://github.com/zhangteng2000/schmeisser-conjecture/blob/9e958aec9d41c387d19a2485fa771895f09eda60/work/uniform_large_degree_constants.py}{\nolinkurl{work/uniform_large_degree_constants.py}}. These comparisons
check the numerical constants used in the large-degree argument
and in Lemma~\ref{lem:local}. Verification requires only the
Python standard library. The launcher refuses to run
under Python's optimization mode, which would disable
assertions. Verification receipts are written to a
temporary directory, so the supplied certificate
files are not altered.

The independent integer verifier checks the finite
interval comparisons required by \eqref{eq:acceptance}.
The reduction from Schmeisser's conjecture to
\eqref{eq:acceptance}, together with the analytic
inequalities used in the reduction, is proved in the main text. The supplement
does not constitute a formal proof-assistant
verification of the complete article.

\end{document}